\documentstyle[12pt, fullpage, amstex]{article}

\newcommand{\Hh}{{\cal H}}

\newcommand{\Tt}{{\cal T}}
\newcommand{\Bb}{{\cal B}}
\newcommand{\Pp}{{\cal P}}

\newcommand{\al}{{\alpha}}

\newcommand{\be}{{\beta}}

\newcommand{\Om}{{\Omega}}
\newcommand{\om}{{\omega}}
\newcommand{\eps}{{\varepsilon}}
\newcommand{\de}{{\delta}}
\newcommand{\De}{{\Delta}}
\newcommand{\ga}{{\gamma}}
\newcommand{\Ga}{{\Gamma}}
\newcommand{\ka}{{\kappa}}
\newcommand{\la}{{\lambda}}

\newcommand{\Diff}{{\rm Diff}}
\newcommand{\Ham}{{\rm Ham}}

\newcommand{\length} {{\rm length}}
\newcommand{\MS}{{\medskip}}

\newcommand{\bl}{{\Big (}}
\newcommand{\br}{{\Big )}}
\def\done{{1\hskip-2.5pt{\rm l}}}
\def\one{{1\hskip-2.5pt{\rm l}}}
\def\NN{{\mathbb N}}
\def\RR{{\mathbb R}}
\def\TT{{\mathbb T}}
\def\ZZ{{\mathbb Z}}
\def\N{{\mathbb N}}
\def\R{{\mathbb R}}
\def\T{{\mathbb T}}
\def\Z{{\mathbb Z}}
\def\bks{{\backslash}}
\def\tet{\theta}		

\def\nek{,\ldots,}
\def\suml{\sum\limits}
\def\prodl{\prod\limits}
\def\intl{\int\limits}
\def\maxl{\max\limits}
\def\tilF{{\widetilde F}}
\def\tilF{{\widetilde F}}
\def\tilp{{\widetilde p}}

\newcommand{\TTd}{\T^d}
\newcommand{\vx}{{ x}}
\newcommand{\va}{{ \alpha}}
\newcommand{\vb}{{ \beta}}
\newcommand{\vh}{{ h}}
\newcommand{\vz}{{ 0}}
\newcommand{\vom}{{ \om}}

\newcommand{\HH}{\mathbf H} 
\newcommand{\im}{\operatorname{Im}} 
\newcommand{\NI}{{\noindent}}

\newcommand{\QED}{\hfill$\Box$\medskip}

\begin{document}
\title{Kick stability in groups and dynamical systems }
\author{
Leonid
Polterovich
\thanks{Supported by THE ISRAEL SCIENCE FOUNDATION
founded by the Israel Academy of Sciences and Humanities.}
\; and Zeev Rudnick
\\ School of Mathematical Sciences\\
Tel Aviv University \\ 69978 Tel Aviv, Israel\\ 
polterov@@math.tau.ac.il \;\;\; rudnick@@math.tau.ac.il\\
}

\date{\today} 

\maketitle

\NI

\begin{abstract} 

We consider a general construction of ``kicked systems'' which extend
the framework of classical dynamics. Let 
$G$ be a group of measure preserving transformations of a 
probability space. Given a one-parameter/cyclic subgroup
(the flow), and any sequence of elements (the kicks) we define
the kicked dynamics on the space by alternately flowing with
given period, then applying a kick. Our main finding is the
following stability phenomenon: the kicked system often
inherits recurrence properties of the  original flow.
We present three main examples. 

1) $G$ is the torus. We show  that for  generic linear flows, and any
sequence of kicks, the trajectories of the kicked system are uniformly
distributed for almost all periods. 

2) $G$ is a discrete subgroup  of $PSL(2,\R)$ acting on the unit tangent
bundle of a Riemann surface. The flow is generated by a single
element of $G$, and we take any bounded sequence of elements of $G$ as
our kicks. We prove that the kicked system is mixing for all
sufficiently large periods if and only if the generator is of
infinite order and is not conjugate to its inverse in $G$.

3) $G$ is the group of Hamiltonian diffeomorphisms of a closed 
symplectic manifold. We assume that the flow is rapidly growing
in the sense of Hofer's norm, and the kicks are bounded. We prove that
for a positive proportion of the periods the kicked system inherits a  
kind of energy conservation law and is thus super-recurrent. 

We use tools of geometric group theory (quasi-morphisms) 
and symplectic topology (Hofer's geometry).

\end{abstract}



\newpage

\section{Introduction and main results}

In the present paper we address the following question:
"How far can a flow be kicked?". 
More precisely,
consider the behavior of a one parameter/cyclic subgroup
of a Lie group under the influence of a sequence of
kicks. The kicks arrive periodically in time. 
The kicks are deterministic, while
the period is chosen "at random". We are interested in
the following stability type question: does the kicked
system inherit some recurrence properties of the 
original one? It turns out that in some  situations
(linear flows on tori, isometries of
$PSL(2,\R)/\Gamma$ , "rapidly growing"
Hamiltonian flows on symplectic manifolds) 
such a stability indeed takes place with positive 
probability even when the kicks are quite large.

\medskip
\noindent
{\it 1.1 Sequential systems}

\medskip
\noindent
Let $G$ be a group. Consider the set 
$G^{\infty} = G \times G \times ...$ of all infinite sequences
${f_*} = \{f_i\}, i \in \N$. Given an action of $G$ on a set
$X$, one can view $f_*$ as a dynamical system (see 
\cite{Bergelson}): The trajectory
of a point $x \in X$ is defined as $\{x_i\}$, where
$x_0 = x,\;x_1 = f_1x,\; x_2 = f_2f_1x$ and so on. We write
$f^{(i)}$ for the
{\it evolution} of $f_*$ given by
 $f_if_{i-1}...f_1$ and set $f^{(0)} = \one$, so
$x_i = f^{(i)}x$.  Note that the constant sequence ${f_*}=\{f\}$
generates iterations of a single map $f$.

Sequences of maps provide a natural framework for the study
of perturbations of the usual settings of dynamical systems 
\footnote{In particular, they are used in various models
of random dynamics \cite{Kifer}, see 1.8 below
for further comments. }. 
We start with the following simple model.
Write $\Tt$ for either $\R$ or $\Z$.
Let 
$(h^t),\; t \in \Tt$ be a one-parameter or cyclic subgroup 
of $G$ which represents
a dynamical system
on $X$. Assume that the system is influenced by a sequence
of kicks $\{\phi_i\} \in G^{\infty}$. The kicks arrive periodically
in time with some positive period $\tau \in \Tt$.
 The kicked
dynamics is described by
a sequential system ${f_*}^{\tau} = \{\phi_i h^{\tau}\}$.
An orbit $\{x_i\}$ of the kicked system looks as follows: 
In order to get $x_i$ from $x_{i-1}$, go with the flow
for the time $\tau$ and then apply the $i$-th kick: $x_i = 
\phi_ih^{\tau}x_{i-1}$. 

One cannot expect general sequential systems to possess
some interesting dynamical properties. For instance if $G$ acts
by measure preserving transformations, a sequential system may
violate the Poincare recurrence theorem etc. However, the kicked
systems described above are very special. 

\medskip
\noindent
{\bf Informal Definition 1.1.A.} A dynamical property 
of a subgroup $(h^t)$ is
called {\it kick stable}, if for every sequence of kicks
$\{\phi_i\}$ (possibly satisfying some mild assumptions)
the kicked system ${f_*}^{\tau}$ inherits this property
for a large set of periods $\tau$.

\medskip
\noindent
In the present paper we discuss various examples and counter-examples
to kick stability. Before presenting a formal definition
(see 1.3.A below) let us specify some dynamical properties we wish
to deal with and consider a number of examples.
The definitions below (possibly with exception of 1.1.C) are
straightforward extensions of standard dynamical notions to the
framework of sequential systems.
Assume that $G$ is a topological group which
acts by measure preserving homeomorphisms
on a topological space $X$ with a Borel probability measure $\mu$.
In 1.1.B - 1.1.D below $X$ is assumed to be compact.
Let ${f_*} \in G^{\infty}$ be a sequential system.

\medskip
\noindent
{\bf 1.1.B. (see \cite{Pjems})} A sequential
system $f_*$ is called {\it strictly ergodic} if
for every continuous function $F$ on $X$ the Birkhoff sums
$${1 \over N} \sum_{i=0}^{N-1} F \circ f^{(i)}$$
converge uniformly to $\int_X F(x)d\mu(x)$.

\medskip
\noindent
{\bf 1.1.C.} A continuous function $F$ is called
 {\it a quasi-integral} of $f_*$ if
$$\lim\sup_{N \to \+\infty}\max_X {1 \over N}\sum_{i=0}^{N-1} F \circ f^{(i)}
> \int_X F(x)d\mu(x).$$

\medskip
\noindent
Clearly these definitions complement each other, that is either $f_*$
is strictly ergodic, or it admits a quasi-integral.
Sometime it is useful to relax the continuity assumption and to work
with characteristic functions $\chi_A$ of measurable 
subsets $A \subset X$.

\medskip
\noindent
{\bf 1.1.D.} If $F = \chi_A$ satisfies 1.1.B we say that $A$ is {\it
ideally recurrent} for $f_*$, and if $F=\chi_A$ satisfies 1.1.C then
$A$ is called {\it super-recurrent} for $f_*$.
Informally speaking,
ideal recurrence means that trajectories of the system
visit $A$ with the frequency $\mu(A)$, while super-recurrence
means that
there exist arbitrarily long finite pieces
of trajectories of $f_*$ which visit $A$ with the frequency $> \mu(A)$.

\medskip
\noindent
{\bf 1.1.E.} A system $f_*$ 
is called  
{\it mixing} if for any two $L^2$-functions
$F$ and $G$ on $X$ the sequence 
$$\int_X F(f^{(i)}x)G(x)d\mu(x)$$
converges 
to $$\int_X F(x)d\mu(x) \int_X G(x) d\mu(x).$$ 

\medskip
\noindent
{\it 1.2 Examples of kick stable flows}

\medskip
\noindent
{\bf Example 1.2.A (stable uniform distribution).} Suppose 
that $G$ is the circle $S^1=\RR/\Z$
which acts on itself by shifts. Consider the flow
$$h^t :S^1 \to S^1,\; h^t x = x+t \mod 1.$$
A characteristic feature of this flow is that for 
irrational values of $\tau$, every trajectory of the sequential
system $\{h^{\tau}\}$ is uniformly distributed in $S^1$, which 
in our language means
that every open interval of $S^1$ is ideally recurrent.
Take an arbitrary sequence of kicks
${\phi_*} \in (S^1)^{\infty}$, and consider the kicked system
${f_*}^{\tau} 
= \{\phi_i h^{\tau}\}$.

\medskip
\noindent
\proclaim Theorem 1.2.B.
There exists a subset
$P \subset (0;+\infty)$ of full Lebesgue measure
such that for every $\tau \in P$ every trajectory of the
kicked system $
{f_*}^{\tau}$ is uniformly distributed in $S^1$. 

\medskip
\noindent
In general, the set $P$ depends on the choice of kicks.
Theorem 1.2.B is a consequence of the Weyl criterion.
We prove it in a more general context of linear flows
on tori in \S 2 below.

\medskip
\noindent
{\bf Example 1.2.C (Stable super-recurrence.)}
Let $(X,\Omega)$ be a closed symplectic manifold.
Denote by $\mu$ the canonical probability measure on $X$.
Let $G=\Ham(X,\Om)$ be the group of all Hamiltonian diffeomorphisms of
$(X,\Omega)$. 
Recall that a symplectic diffeomorphism is called Hamiltonian if
it can be included into a time-dependent Hamiltonian flow.
A Hamiltonian function $H: X \times [0;1] \to \R$ is called
{\bf normalized} if
$\int_X H(x,t) d\mu(x) = 0$ for all $t \in [0;1]$. Every Hamiltonian
diffeomorphism $h \in G$ can be 
written as $h = h^1$, where $h^t$ is the Hamiltonian
flow generated by some normalized Hamiltonian 
$H$. In this case we say that $h$ is generated
by $H$. Define a function ${\bar \rho} : G \to [0;+\infty)$ as follows:
$${\bar \rho} (h)  = \inf \int_0^1 \max_{x \in X} H(x,t) - 
\min_{x \in X} H(x,t) dt,$$
where the infimum is taken over all normalized Hamiltonians 
$H $ which generate $h$.
The function $\bar \rho$ is known as 
{\it Hofer's norm}
(see \cite{Hofer},\cite{Picm},\cite{Pbook}).
Let us introduce the following important notion.
A sequence $\{\phi_i\} \in G^{\infty}$ is called {\bf bounded}
if the sequence $\{{\bar \rho}(\phi_i)\}$
is bounded.
 
We illustrate stable super-recurrence 
in the simplest case 
when 
$X = S^2$ is the unit sphere in $\R^3$ endowed with the
induced Euclidean area form. In this case every diffeomorphism
which preserves $\mu$   
and the orientation is Hamiltonian.
Every one-parameter subgroup $(h^t)$
of $G$ is simply an autonomous Hamiltonian flow generated
by a (uniquely defined) time-independent Hamiltonian function
$H \in \Hh$. Clearly the function $H$ is an integral of motion
(the energy conservation law!), and in particular
each subset
$$A_{\epsilon} = \{ H > (1-\epsilon)\max H \},\; \epsilon \in (0;1)$$ 
is invariant under the flow,
and thus super-recurrent in the sense of Definition 1.1.D.
Take an arbitrary bounded sequence 
of kicks 
$\{\phi_i\} \in G^{\infty}$.
\footnote
{
 For instance, one can think
that all $\phi_i$'s are conjugate to elements of some compact
(in the $C^{\infty}$-topology)
subset of $G$.
}
Consider the kicked system $
{f_*}^{\tau} 
= 
\{\phi_i h^{\tau}\}$.
Fix $\epsilon \in (0;1)$. 

\medskip
\noindent
\proclaim Theorem 1.2.D. 
Suppose that the Hamiltonian $H$ is non-constant, and its maximum
set 
$\{H = \max H\}$ contains a simple closed curve which divides
the sphere into two discs of equal areas. Then
there exists a subset $P \subset (0;+\infty)$
whose complement has a finite measure and such that for every
$\tau \in P$ 
\begin{itemize}
\item{}the Hamiltonian $H$ is a quasi-integral
of the kicked system ${f_*}^{\tau}$;
\item{}the set $A_{\epsilon}$ is super-recurrent for 
${f_*}^{\tau}$.
\end{itemize}

\medskip
\noindent
For instance, 
consider the flow $(h^t)$ which 
rotates each point $(x,y,z) \in S^2$ with the velocity $z$ around
the $z$-axis. 
In Euclidean coordinates $(x,y,z)$
on $\R^3$ it is given by
$$h^t (x,y,z) = (x \cos (2\pi tz) - y \sin (2\pi tz), 
x \sin (2\pi tz) + 
y \cos (2\pi tz),z). \leqno(1.2.E) $$

It is generated by the Hamiltonian
function $H(x,y,z) = -z^2+ {1\over 3}$ 
whose maximum set coincides with the
equator $\{z=0\}$. Therefore the theorem above is applicable,
and in particular each annulus $\{|z| < \varepsilon\}$
is super-recurrent for the kicked system $f_*^\tau$ for all values
of the period $\tau \in P$.

\medskip
\noindent
 Theorem 1.2.D is proved in 4.3 below.
Even in this simple case the proof we have is
not elementary - it is
based on 
powerful methods of modern symplectic topology.
\footnote
{Dima Burago pointed out that in the case when the sequence
of kicks $\{\phi_i\}$ is contained in a compact subset of $G$,
the conclusions of 1.2.D can be checked by soft methods
(compare with the previous footnote).}
We refer to section 1.6 for
various extensions of this theorem to more general
symplectic manifolds.   

\medskip
\noindent
{\bf Example 1.2.F. (stable mixing)}
Let $\Gamma \subset PSL(2,\R)$ be a lattice, 
that is a discrete subgroup such that
the Haar measure $\mu$ of
the quotient space $X = PSL (2,\R) /\Gamma$ is finite.
Let $G \subset PSL(2,\R)$ be a {\em discrete} subgroup, 
and consider the left action of $G$ on $X$.
Let $h \in G$ be an element of infinite order. The Howe-Moore theorem
\cite{Zimmer}
(see also 1.5.A below) yields that $h$ is a mixing transformation
of $X$. 
Let $\phi_* = \{\phi_i\}$ be an arbitrary sequence from $G^{\infty}$ 
which represents a finite number of conjugacy classes in $G$.
\footnote
{For instance, if $\Gamma = PSL (2,\Z)$ this assumption
holds when all $\phi_i$'s are non-parabolic and their
traces are bounded.}
Consider the kicked system ${f_*}^\tau = \{\phi_i h^{\tau}\}$, 
where $\tau \in \N$.
We say that $h$ is {\it stably mixing} if for every 
sequence $\phi_*$ as above
there exists $\tau_0 > 0$ 
such that the kicked system 
${f_*}^\tau$ is mixing for all $\tau > \tau_0$.
The next result gives a complete description of stably mixing elements
of $G$ in purely algebraic terms.

\medskip
\noindent
\proclaim Theorem 1.2.G. Let $h \in G$ be an element of infinite
order. The following conditions are equivalent:
\begin{itemize}
\item[(i)] $h$ is stably mixing on $X$; 
\item[(ii)] $h$ is not conjugate
to its inverse $h^{-1}$ in $G$.
\end{itemize}

\medskip
\noindent
For instance, it is easy to see  that the time one map of the 
horocycle flow on $X$ which is given by the matrix
$$ 
\left (
\begin{array}{cc}
1& 1\\
0& 1\\
\end{array}
\right ) 
$$
is not conjugate to its inverse already in $PSL(2,\R)$ 
and thus is stably mixing. 
On the other hand in 
$G = PSL(2,\Z)$ every symmetric matrix
is conjugate to its inverse by the involution
$$
\left (
\begin{array}{cc}
0& -1\\
1& 0\\
\end{array}
\right )
,
$$
and thus is not stably mixing. 
A complete description of $PSL(2,\Z)$-matrices
which are conjugate to their inverses is unknown
(cf. \cite{Baake}). 
We refer to 1.7 for the
fairly general discussion on the
effect of time-reversing symmetry on the kick stability.
The proof of Theorem 1.2.G is given in  1.5, 1.7 and \S 3 below.
\footnote{
The condition which tells that the sequence
of kicks $\{\phi_i\}$ represents a finite number of conjugacy classes
in $G$ is essential as the next simple example shows.
Let $G=PSL(2,\Z)$ and $h=\left(\smallmatrix
1&1\\0&1\endsmallmatrix\right)$.  
Choose a surjective
function $\alpha : \N \to \N$ such that  $\alpha^{-1}(\tau)$
is an infinite subset for every $\tau \in \N$. 
Put $q_k = k\alpha(k) -(k-1)\alpha(k-1),$
where $k \in \N$. Consider a sequence of kicks
$\phi_k = h^{-q_k}$. Clearly it represents 
an infinite number of different
conjugacy classes in $PSL(2,\Z)$.
The evolution of the kicked system
is given by
$f^{(k)}(\tau) =h^{ k(\tau-\alpha (k))}.$ Thus for every $\tau \in \N$
equality $f^{(k)}(\tau) = \one$ holds for an infinite number of $k$'s,
so the kicked system is not mixing. 
}

\vfill\eject

\medskip
\noindent
{\it 1.3 Kick stability.}

\medskip
\noindent
Now we are ready to give a formal definition of kick stability
which provides a unified framework to examples considered above.
In what follows $\Tt$ stands either for $\R$ or for $\Z$,
and $\Tt_+ = \{t \in \Tt \;{\big |}\; t > 0\}$. We fix a class
$\Bb$ of subsets of $\Tt_+$ of "large measure". 
For instance
we can assume that when $\Tt = \R$ (resp. $\Tt = \Z$) the
class $\Bb$ consists of all subsets of $\Tt_+$ whose complement
has finite Lebesgue measure (resp. is a finite subset). 

Let $G$ be a group, and let $(h^t),\; t \in \Tt$ be a subgroup
which is assumed
to be either one parameter ($\Tt = \R$) or cyclic
($\Tt = \Z$).
This subgroup represents the unperturbed dynamical system with 
continuous or discrete time. 
Fix a subset $\Pp \subset G^{\infty}$
which should be thought of as the set of all sequential systems
with a given property (P). 
We say that our subgroup $(h^t)$ has property (P) if
the set
$$\{\tau \in \Tt_+ \;{\big |}\; {\rm the \; sequence}\; \{h^{\tau}\} \in 
\Pp\}$$
has large measure, that is belongs to $\Bb$. 

Let $\Phi \subset G^{\infty}$ be a set of admissible
kicks. Take ${\phi_*} \in \Phi$ and
write ${f_*}^\tau ({\phi_*})$ for the kicked
system $\{\phi_i h^{\tau}\}$. Denote by $P({\phi_*})$
the set of all positive values of the period $\tau$ such that
the kicked system has property (P):
$$P({\phi_*}) = \{\tau \in \Tt_+ \; {\big |}\;
{f_*}^{\tau} ({\phi_*}) \in \Pp \}.$$

\medskip
\noindent
{\bf Formal Definition 1.3.A.} The property (P) of
subgroup $(h^t)$ is {\it kick stable} if for every
admissible sequence of kicks ${\phi_*} \in \Phi$
the set $P({\phi_*})$ of "good" periods belongs
to the class $\Bb$ of subsets of large measure.

\medskip
\noindent

Let us illustrate this definition. In Example 1.2.A above,
all kicks are admissible so
the set $\Phi$ coincides with $(S^1)^{\infty}$. The set
$\Pp$ consists of all sequences from $(S^1)^{\infty}$ whose orbits
are uniformly distributed in $S^1$. Theorem 1.2.B implies
that the property {\it all orbits are uniformly
distributed in $S^1$} is kick stable for $(h^t)$. In Example 1.2.C
the set $\Phi$ of admissible kicks consists of all sequences which
are bounded in Hofer's norm.
The set $\Pp$ is formed by all sequential
systems for which the set $A_{\epsilon}$ is super-recurrent.
Theorem 1.2.D states that the property {\it  $A_{\epsilon}$ is
 a super-recurrent set} is kick stable for $(h^t)$. In Example 1.2.F
all sequences of kicks which represent a finite number of conjugacy
classes
 are admissible, and the set $\Pp$ consists of all mixing
systems. Theorem 1.2.G tells us when {\it  mixing} is a kick stable
property of the cyclic subgroup $(h^t)$.

\medskip
\noindent
{\it 1.4 Sub-additive functions}

\medskip
\noindent
We do not know of an argument
which provides a unified explanation of the kick-stability 
phenomenon in all the examples
presented in section 1.2 above. Interestingly enough, however
 that our approaches to
Theorems 1.2.D and 1.2.G have a common ingredient. Namely,
the desired kick
stability is closely related to the geometric behaviour 
of the corresponding
subgroups at infinity. 

\medskip
\noindent
{\bf Definition 1.4.A.} Let $G$ be a group. A function
$\rho : G \to [0;+\infty)$ is called {\it sub-additive}
if there exists a number $C \geq 0$ such that
the following holds:
\begin{itemize}
\item{} $|\rho(hgh^{-1}) - \rho(g)| \leq C$ for all $g,h \in G$;
\item{} $\rho(gh) \leq \rho(g)+\rho(h) + C$ for all $g,h \in G$.
\end{itemize}

\medskip
\noindent
If $\rho$ is sub-additive then, as is well known,
for every $h \in G$, the
limit $$ \rho_{\infty} (h) = \lim_{n \to +\infty}{1 \over n} 
\rho_{\infty}(h^n)$$
does exist. In the next two sections we 
discuss two applications of sub-additive functions
on groups to kick stability.

\medskip
\noindent
{\it 1.5 Stable mixing and quasi-morphisms}

\medskip
\noindent
In this section we outline the proof of the sufficient condition
of stable mixing for  an element of the discrete group $G$
given in 1.2.G above.

We start with a following more general situation.
Suppose that $D$ is a non-compact simple Lie group with finite center.
Let $\Gamma \subset D$ be a lattice, that is a discrete subgroup
such that the Haar measure of the quotient space $X = D/ \Gamma$
is finite. 
The group $D$ 
acts on $X$ on the left by transformations
preserving the Haar measure. The key ingredient 
of our approach to stable mixing is the 
Howe-Moore theorem which provides a link between geometry
and dynamics of sequential systems in this setting. In order to
formulate it we need the following notions.
Let $\{f^{(i)}\}$ be a sequence
of elements of $D$. We say that $\{f^{(i)}\}$ {\it goes to infinity}
if for every compact subset $Q \subset D$ there exists
$i_0$ such that $f^{(i)} \notin   Q$ for all $i > i_0$. From general
considerations, it follows that a mixing sequence of elements of $D$
necessarily goes to infinity~\footnote{ 
Indeed, otherwise there is a subsequence which is
{\em bounded} and therefore there is a {\em convergent} subsequence,
which is also mixing. Thus (after passing to this subsequence) 
we have a sequence $f^{(i)}\to f_\infty$  and that is still
mixing. Now take a real valued
function $G$, not identically zero, with $\int_X G(x)d\mu(x)=0$ and set 
$F(x):=G(f_\infty^{-1} x)$. By the mixing property of the sequence, we
have 
$$
\int_X F(f^{(i)}x)G(x)d\mu(x) \to \int_X F(x)d\mu(x)\int_X G(x)d\mu(x)= 0
$$
while because $f^{(i)}\to f_\infty$, we have 
$$
\int_X F(f^{(i)}x)G(x)d\mu(x) \to \int_X F(f_\infty x)G(x)d\mu(x)
=\int_X G(x)^2d\mu(x) >0
$$
which gives a contradiction.}. 
Conversely, we have:   

\medskip
\noindent
\proclaim Howe-Moore Theorem 1.5.A. (\cite{Zimmer}). Let $f_*$ be a sequential
system from $D^{\infty}$. If
its evolution $\{f^{(i)}\}$ goes to infinity 
then $f_*$ is mixing.

\medskip
\noindent
Let $G \subset D$ is a discrete group. Below we focus on the action
of $G$ on $X$.
 Consider the cyclic subgroup
generated by an element 
$h \in G$. Let $\phi_*$ be an arbitrary sequence from $G^{\infty}$ whose
entries represent a finite number of conjugacy
classes in $G$. Consider the kicked system
${f_*}^\tau = \{\phi_i h^\tau\}$, where $\tau \in \N$.
The next result provides a sufficient 
condition for kick stability of mixing
for our subgroup.

\medskip
\noindent
\proclaim Theorem 1.5.B. Assume that there exists a sub-additive function
$\rho$ on $G$
such that $\rho_{\infty}(h) > 0$. Then there exists $\tau_0 >0$
such that the kicked system ${f_*}^\tau$ 
is mixing for all $\tau > \tau_0$.

\medskip
\noindent
Here $\tau_0$ depends on $\phi_*$. 
It turns out that the geometric assumption $\rho_{\infty}(h) > 0$ 
guarantees that
for large periods $\tau$
the evolution of the kicked system 
goes to infinity , thus the statement follows
from 1.5.A. 
The details of this argument are given in 3.1 below.

There exists a useful class of sub-additive functions which
arise naturally in the bounded cohomology theory of discrete
groups (see \cite{Brooks},\cite{Barge},\cite{Picaud}).

\medskip
\noindent
{\bf Definition 1.5.C.} A function $r: G \to \R$ is called
a {\it quasi-morphism} if there exists a constant $C > 0$
such that
$$|r(gh)-r(g)-r(h)| \leq C$$
for all $g,h \in G$. 

\medskip
\noindent
Given a quasi-morphism $r$ and an element $g \in G$ ,
there exists the limit
$$r_{\infty}(g) = \lim_{n\to +\infty} {r(g^n)\over n}.$$
Note that $r_\infty$ is homogeneous, that is
$r_\infty(g^k)=kr_\infty(g)$. Moreover, if $g$ has finite order then
$r_\infty(g)=0$. 
It follows immediately from the definition that if $r$ is a quasi-morphism
then the function $\rho(g) = |r(g)|$ is sub-additive, and moreover
$\rho_{\infty}(g) = |r_{\infty}(g)|$ for all $g \in G$. 

\medskip
\noindent
\proclaim Theorem 1.5.D. Let $G\subset PSL(2,\R)$ be a discrete group
and $h$ an element of infinite
order in $G$. The following conditions are equivalent:
\begin{itemize}
\item[(i)] there exists a quasi-morphism $r:G \to \R$ such that
$r_{\infty}(h) > 0$;
\item[(ii)] $h$ is not conjugate to its inverse $h^{-1}$ in $G$.
\end{itemize}

\medskip
\noindent
{\bf Proof of ``1.2.G(ii) implies 1.2.G(i)'':}
The desired statement is an immediate consequence
of 1.5.D and 1.5.B.
\QED

\medskip
\noindent
Theorem 1.5.D is proved in 3.2 below 
(see Remark 3.2.F for references and generalizations
of this result).

\medskip
\noindent
{\it 1.6 Stable super-recurrence in Hamiltonian dynamics}

\medskip
\noindent
Let $(X,\Omega)$ be a closed symplectic manifold.
Denote by $\mu$ the canonical probability measure on $X$.
Let $G=\Ham(X,\Om)$ be the group of all Hamiltonian diffeomorphisms of
$(X,\Omega)$. Define {\it the positive part of Hofer's norm}
$\rho : G \to [0;+\infty)$ as follows:
$$\rho (h) = \inf \int_0^1 \max_{x \in X} H(x,t) dt,$$
where the infimum is taken over all normalized Hamiltonian
functions $H : X \times [0;1] \to \R$ which generate $h$
(cf. 1.2.C above).
It is an easy exercise to check that $\rho$ is sub-additive
(here the constant $C$ of 1.4.A is simply $0$).
Moreover, the following obvious inequality
holds:${\bar \rho} (h) \geq \rho(h) + \rho (h^{-1})$
for all $h \in G$.
\footnote
{In fact, in all known examples one has the equality!}

Let $(h^t)$ be a one-parameter subgroup of $G$. It is generated by
some uniquely defined time-independent 
Hamiltonian $H : X \to \R$ with zero
mean.
The law of energy conservation yields that $H$ is an integral of motion.
Note that 
$\rho (h^t) \leq t\max H$ for all $t > 0$, and thus we have
$\rho_{\infty}(h^1) \leq \max H$.
Let $\phi_* \in G^{\infty}$  be an arbitrary bounded sequence 
of kicks (see 1.2.C above).
Consider the kicked system
${f_*}^\tau = \{\phi_i g^\tau\}.$

\medskip
\noindent
\proclaim Theorem 1.6.A. (stable energy conservation law). 
Assume that $\rho_{\infty}(h^1) > 0$.
There exists a subset $P \subset (0;+\infty)$ 
of 
density at least 
$${\rho_{\infty}(h^1)}\over {\max H}$$
such that for every $\tau \in P$ the Hamiltonian
$H$ is a quasi-integral of the kicked system
${f_*}^\tau$.

\medskip
\noindent
This result gives a probabilistic interpretation
for a purely geometric quantity $\rho_{\infty}(h^1)$.
Further, this quantity contains 
interesting information about kick stable super-recurrence of 
$(h^t)$-invariant sets of the form 
$$A_{\epsilon} = \{H > \epsilon \max H\}.$$
We present here two sample results. 

\medskip
\noindent
\proclaim Theorem 1.6.B. Assume that $\rho_{\infty} (h^1) = \max H$.
Fix $\epsilon > 0$. There exists a subset $P \subset (0;+\infty)$
of density 1 such that 
for every $\tau \in P$
the set $A_{\epsilon}$ is super-recurrent for 
the kicked system ${f_*}^\tau$.

\medskip
\noindent
One can check (see 4.2, 4.3 below) that in Example 1.2.C above
the assumption 
$\rho_{\infty} (h^1) = \max H$
is satisfied, thus 
(a slightly weaker version of)
Theorem 1.2.D follows from 1.6.B.

\medskip
\noindent
\proclaim Theorem 1.6.C. Assume that 
$\rho_{\infty} (h^1) \geq 0.9 \max H$,
and $\max H = -\min H$. Then there exists a 
subset $P \subset (0;+\infty)$
of density at least 0.4 such that for every 
$\tau \in P$ the set $A_{0.4}$
is super-recurrent for the kicked system ${f_*}^\tau$.

\medskip
\noindent
We prove more general versions of these theorems in \S 4 below. 

\medskip
\noindent
{\it 1.7 An obstruction to kick stability}

\medskip
\noindent
We describe here a method of constructing
kick unstable systems in a number of interesting situations.
Let $(h^t), \; t \in \Tt$ be a one-parameter/cyclic
subgroup of a group $G$.

\medskip
\noindent
{\bf Definition 1.7.A.}  An element $\theta \in G$ is called
{\it a time-reversing symmetry}
\footnote{See \cite{Lamb} for a discussion on time-reversing
symmetries and their impact on dynamics.}
for $(h^t)$ if 
$\theta h^t \theta^{-1} = h^{-t}$ for all $t \in \Tt$.
\medskip
\noindent
Let us introduce the following notation.
Given a system $f_* \in G^{\infty}$ we write $f^* = \{f^{(i)}\}$
for its evolution
$f^{(i)} = f_i...f_1$. 

\medskip
\noindent
{\bf 1.7.B. Creating periodic behaviour.} 
Suppose that $(h^t)$ admits
a time-reversing symmetry $\theta$.
Take a sequence of kicks 
$$\phi_* = \{\theta^{-1},\theta,\theta^{-1},\theta,...\}.$$
Consider the
kicked system ${f_*}(\tau) = \{\phi_i h^\tau\}$.
Its evolution $f^*(\tau)$ is 2-periodic:
$$ f^*(\tau) = \{\theta^{-1}h^{\tau},\one,
\theta^{-1}h^\tau,\one,...\}.$$
Assume that $G$ acts by measure-preserving homeomorphisms
on a topological space $X$. 
Obviously, for all $\tau \in \Tt$
every subset $A$ of $X$ with $\mu (A) < 0.5$
is super-recurrent for the kicked system $f_*(\tau)$.

\medskip
\noindent
{\bf Example 1.7.C. (cf. 1.2.A above).}
Consider the action of the orthogonal
group $O(2)$ on the circle $S^1$. The uniform distribution 
property of the flow $h^t x = x +t$ is not kick stable
for kicks from $(O(2))^{\infty}$.
Indeed, the transformation $\theta: x \to -x$ is a
time-reversing symmetry for $h^t$.  Comparing this with Theorem
1.2.B above we see that the flow $h^t$ loses stability when
one replaces the group $S^1 = SO(2)$ by a larger group
$O(2)$. This stability breaking mechanism can be observed 
in many other situations (see 1.7.H and 4.6.A below).

\medskip
\noindent
{\bf 1.7.D.} Suppose now in addition that every point $x \in X$ admits
a nested system of open neighborhoods 
$U_{\delta},\; \delta \in (0;\delta_0)$ such that
$ \cap_{\delta} U_{\delta} = \{x\}$ and $\mu(U_{\delta}) = \delta$.
We claim that for every $\tau$
the kicked system $f_*(\tau)$ constructed above
is not mixing.
Indeed, fix arbitrary $\tau \in \Tt$, and take $\delta > 0$ small
enough.
Since  $\theta \neq h^{\tau}$
there exists an open subset
$U$ of $X$ of measure $\delta$ such that
$\theta^{-1}h^\tau U \cap U = \emptyset$.
Let $F$ be a characteristic function of $U$.
The only limit points of the sequence
$\int_X F(f^{(i)}(\tau)x)F(x)d\mu(x)$ are
$0$ and $\delta$, while
${\big (}\int_X F(x)d\mu(x) {\big )}^2 = \delta^2$.
Therefore the kicked system is not mixing.

\medskip
\noindent
{\bf Proof of ``1.2.G(i) implies 1.2.G(ii)'':}
Consider the action of $G$ on $ X= PSL(2,\R)/\Gamma$.
Let $h \in G$ be an element of infinite order.
If $h$ is conjugate to its inverse (that is it admits a
time-reversing symmetry) then the argument 1.7.D above shows that
$h$ is not stably mixing. This completes the proof of 1.2.G.
\QED 

\medskip
\noindent
{\bf 1.7.E. Creating random behaviour.}
Exactly as in examples above, a time-reversing symmetry
provides an obstruction to stable super-recurrence in Hamiltonian
dynamics. We work in the setting of 1.6, assuming in addition
that the group of Hamiltonian diffeomorphisms
 $G$ is $C^{\infty}$-closed in 
the group of all smooth
diffeomorphisms
of $X$. 
\footnote{
This is true for many symplectic manifolds, 
for instance if $H^1(X,\R) = 0$
or if the cohomology class of the symplectic form $\Omega$ is 
rational. The
famous Flux Conjecture (see e.g. \cite{Lalonde}) which states that 
this is true for all
$(X,\Omega)$, is still open.
}
The next result
is proved in 4.6 below.

\medskip
\noindent
\proclaim Proposition 1.7.F. Assume that $(h^t)$ is a one-parameter
subgroup of the group $G$ of Hamiltonian diffeomorphisms of $X$
which admits a time-reversing symmetry. Then there exists
a bounded sequence $\phi_* =\{\phi_i\} \in G^{\infty}$ 
such that 
the kicked system
$\{\phi_i h^\tau\}$ is strictly ergodic for all $\tau \in (0;+\infty)$.

\medskip
\noindent
{\bf Example 1.7.G (cf. 1.2.C above.)} Let $X$ be the 2-sphere,
and let $h^t$ be the circle action which rotates the sphere
with constant speed around the $z$-axes. This flow has no
stably super-recurrent invariant sets with piece-wise
smooth boundary. Indeed, it has a time- reversing symmetry
(for instance, the reflection around the $x$-axes).

\medskip
\noindent
{\bf Example 1.7.H (cf. 1.2.C above.)}
Assume now that $h^t$ is the flow given by (1.2.E) above.
Instead of the group of all Hamiltonian diffeomorphisms of 
$S^2$ consider the larger group of all measure-preserving
diffeomorphisms
\footnote
{Besides the Hamiltonian ones,
this group includes orientation reversing diffeomorphisms
preserving the measure.}. 
It turns out that after such an enlarging of the ambient
 group, the flow $(h^t)$ gets a time-reversing symmetry-
the reflection over the $(x,y)$-plane. A little modification
of 1.7.F above yields that $(h^t)$ loses stable super-recurrence
(cf. 1.7.C above).

\medskip
\noindent
{\it 1.8. Discussion and open problems}

\medskip
\noindent
{\it 1.8.A. Historical and bibliographical remarks.}

\medskip
\noindent
{\bf Sequential systems} arise naturally in random dynamics
\cite{Kifer}. In the deterministic language they form a particular
case of a skew-product. The study of ergodic properties of
individual sequential systems was initiated by Bergelson and Berend 
\cite{Bergelson}. Some special classes of sequential systems were known
for a long time. For instance, sequences of M\"obius transformations
were considered in connection with analytic continued fractions
as well as with the discrete Schr\"odinger operator. Very recently
Zeghib \cite{Zeghib} investigated sequences of isometries
of a Lorentz manifold. He found that the corresponding dynamics
is closely related to asymptotic geometry of the isometries
group (cf. 1.5-1.6 above). It is clear that 
sequential dynamics have not been 
systematically studied yet, and many natural and interesting
questions are still unexplored.

{\bf Kicked systems} were intensively studied by physicists in
the classical (= non-sequential) framework. This
class of systems includes a number of famous maps which
attracted a lot of attention in conservative
chaotic dynamics, such as the kicked top
\cite{Haake}, the kicked harmonic oscillator (or the
Henon map \cite{Heagy}) and the kicked rotator 
(or the standard map \cite{Haari}).
Physicists, however considered these maps from a viewpoint which
essentially differs from ours. 
In order to illustrate the difference, let us return to the
 the flow $(h^t)$ given by (1.2.E)
above. Take the constant sequence of kicks
$\phi_i \equiv \phi$, where $\phi (x,y,z) = (-z,y,x)$.
The corresponding kicked system describes iterations
of the single map $\phi h^\tau$, which is nothing else but
{\it the kicked top
map} (see \cite{Haake}).  
In contrast to our setting, in the physics 
literature $\phi$ describes the top,
while $h^\tau$ stands for the kick! Further, 
computer experiments performed in \cite{Haake} suggest that for large
values of $\tau$ the kicked system is chaotic: 
a generic trajectory of the kicked top map
is uniformly distributed on a huge subset of the sphere.
On the other hand, Theorem 1.2.D above guarantees super-recurrence
for all positive $\tau$ outside the set of finite measure.
We arrive at a seemingly paradoxical situation: chaos coexists
with super-recurrence. The resolution of this paradox is as follows.
The deterministic
behaviour of a {\it generic} trajectory of the unperturbed system
is a kick unstable property. This kick instability, reflected
in the transition to chaos for large values of $\tau$, is the
main attraction for physicists. The main message of our theory
is that even for large $\tau$ some (non-generic!) trajectories
are still super-recurrent.
Indeed, our
notion of super-recurrence takes into account behaviour of
{\it all} (vs. almost-all) trajectories of a system.
For instance, the original kicked top map
is super-recurrent simply due to existence
a fixed point $(1,1,0)$. It surprising however that 
this super-recurrence persists when the constant sequence
of kicks is replaced by an arbitrary
bounded non-constant sequence. 
Let us mention also that this conflict between
"all" and "almost all" is a reflection of the striking
difference between $L^p-$ and $C^0-$ measurements in symplectic
topology, see for instance \cite{Eliashberg},\cite{Pbook}. 
A similar analysis of the standard map
(which we omit here)  also leads to some "paradoxical" conclusions
which are far from being understood yet. 

{\bf Kick stability}, the central notion promoted
in the present paper, lies in a long series of attempts
to formalize robustness of dynamical and ergodic properties
of flows and maps (such as recent works on stable ergodicity,
see \cite{Field} and references therein).
 Among its cousins one may recall  stochastic stability
\cite{Viana} which naively speaking means stability with respect
to small random sequential kicks. It would be interesting to
compare kick stability (where the kicks are deterministic
and not assumed to be small) with  stochastic stability
in more detail.

\medskip
\noindent
{\it 1.8.B. Discrete vs. continuous} 

\medskip
\noindent
In basic examples considered above- 
uniform distribution on $S^1$,
super-recurrence in Hamiltonian dynamics and
mixing on $PSL(2,\R)/\Gamma$ - one can address
both discrete and continuous versions of
the question on kick stability. In the present
paper we worked out the continuous versions
in the first two examples, and the discrete version
in the last one. What happens with the remaining cases?
It turns out that cyclic subgroups of $S^1$ do not
have kick stable uniform distribution property.
This was noticed by Dima Burago whose argument is presented
in 2.3 below.  Further, nothing is known to us
about kick stable super-recurrence for cyclic subgroups of
the group of Hamiltonian diffeomorphisms of a symplectic
manifold. It would be interesting to 
make some progress in this direction.
Finally, we arrive at the  kick stability question
for flows
on
$PSL(2,\R)/\Gamma$. It deserves a special discussion.

Consider
the action of $G=PSL(2,\R)$ on $X = PSL(2,\R)/\Gamma$,
where $\Gamma$ is a lattice. It follows from 1.5.A that
every non-compact one-parameter subgroup of $G$
is mixing. Is this property kick stable? For instance, 
the geodesic flow 
$$h^t = 
\left (
\begin{array}{cc}
e^t& 0\\
0& e^{-t}\\
\end{array}
\right )
$$
is not stably mixing. Indeed, it is given by a symmetric
matrix and thus admits a time-reversing symmetry (see discussion
in 1.7.D above). 

{\medskip
\noindent
{\bf Question 1.8.C.} Is the horocycle flow
$$h^t = 
\left (
\begin{array}{cc}
1& t\\
0& 1\\
\end{array}
\right )
$$
stably mixing?

\medskip
\noindent
Note that the horocycle flow does not admit a time reversing
symmetry in $PSL(2,\R)$. Being unable to find the complete
answer, we present some partial results (suggesting
the affirmative solution) and more discussion in the Appendix to
\S 3 below.
In particular, we present a link between this problem
and spectral theory for the discrete Shr\"odinger equation.

\medskip
\noindent
{\it 1.8.D. Generalizations to other dynamical systems.}

\medskip
\noindent
It would be interesting to investigate stable mixing
for lattices in semi-simple Lie groups of higher
rank. Does there exists a complete description of stably
mixing elements (if any) similar to Theorem 1.2.G above?
Recently Burger and Monod \cite{Burger-Monod1, Burger-Monod2} showed
that unlike the rank-one case, there are few quasi-morphisms of
lattices in higher rank groups. Thus new ideas are needed in order to
understand that case.


Another potential source of kick-stable systems
might be provided by hyperbolic theory.  
Here is a warm up question: under which conditions
is the existence of a hyperbolic attractor  a kick stable
property? 

\vfill\eject

\bigskip
{\bf

\centerline{CONTENTS}

\medskip

Section 2. Linear flows on tori.

\medskip

Section 3. Detecting stable mixing on $PSL(2,\R) /\Gamma$.

\medskip

Section 4. Stable super-recurrence in Hamiltonian dynamics.}

\medskip

\noindent{\bf Acknowledgments:} We thank Dima Burago
for a number of important critical remarks and suggestions,
and in particular for explaining us a construction
presented in 2.3 below. We are grateful to Anna Dioubina
for pointing out an error in the original version of
Remark 3.3.E, as well as to Brian Bowditch and Iosif Polterovich
for illuminating consultations on reference [EF].
We thank
Leonid Pastur and Misha Sodin
for various useful discussions. 
\vfill\eject

\section{Linear flows on tori} 

\medskip
\noindent
{\it 2.1. Generic flows are kick stable.}

\medskip
\noindent
In this section we consider kick-stability for linear flows on the
$d$-dimensional torus $\TTd=\R^d/\Z^d$. It turns out that some very
classical results on the ``metric'' theory of uniform distribution
immediately imply that the kicked flows are stably uniformly
distributed for almost all periods. 

Precisely, for $\vom\in \R^d$ consider the  ``Kronecker map'' on
the torus $\TTd$ given by $\vx\mapsto \vx+\vom \mod 1$. It gives a
one-parameter group ($t \in \R$)
$$
h^t(\vx)=\vx+t\vom \mod 1 \;.
$$
Given $\vb_i\in \TTd$, we get ``kicks'' $\phi_i(\vx)=\vx+\vb_i \mod 1$. 
 Put $f_i^{\tau}=\phi_i  h^{\tau}$, so that
$f_i^{\tau}(\vx)=\vx+\vb_i+\tau\vom \mod 1$. 
The evolution of the kicked system is then given by 
$$
f^{(k)} = f_k^{\tau} f_{k-1}^{\tau} \dots 
f_1^{\tau}:\vx \mapsto  \vx +\va_k+k\tau\vom \mod 1
$$ 
where $\va_k=\vb_1+\dots +\vb_k$. 

We will say that $\vom=(\om_1,\dots, \om_d)$ is 
a ``generic'' vector if its 
components $\om_i$ are linearly independent over the rationals.

\medskip
\noindent
\proclaim Theorem 2.1.A.
Suppose that $\vom$ is a ``generic'' vector. Then for
almost all periods $\tau$, the orbits  
$\{f^{(k)}(\vx)\}_{k=1}^\infty$ are
uniformly distributed in $\TTd$.

\medskip
\noindent
This result is a version of a ``metric theorem'' of Weyl from
1916 (see \cite{KN}).
 For the sake of completeness
we recall the argument. In the following,
given two sequences $a(N)$ and $b(N)$
we will use the notation 
$a(N)\ll b(N)$  to mean that there is a constant
$c>0$ so that $a(N)\leq cb(N)$ for all $N$ sufficiently large.

\medskip
\noindent
{\it 2.2. Proof of 2.1.A:}

\medskip
\noindent
Since $\TTd$ acts transitively on itself by translations, it suffices
to consider the base point $\vz$ instead of arbitrary $\vx$. 
It also suffices to fix a finite interval $[a,b]$ and show that the
result holds for almost all $\tau\in[a,b]$.

Define normalized ``Weyl sums'' 
$$
 S_{\vh}(N,\tau):= \frac 1N\sum_{k=1}^N 
e^{2\pi i\vh \cdot f^{(k)}(\vz)}
$$
The basic tool is Weyl's observation that for uniform distribution, it
suffices to show that the normalized Weyl sums converge to zero for
all integer vectors $\vh\neq \vz$. 

To do that, one shows (see below) that for fixed $\vh\neq \vz$, one has 
$$ 
\int_a^b \left| S_{\vh}(N,\tau) \right|^2 d\tau \ll\frac {\log N}N.  
\;\leqno(2.2.A)
$$
Thus for the sequence of  squares $N=n^2$, we have 
$$
\sum_{n=1}^\infty \int_a^b \left| S_{\vh}(n^2,\tau) \right|^2 d\tau
< \infty \; .
$$
By Fatou's lemma, it follows that 
$\sum_{n=1}^\infty \left| S_{\vh}(n^2,\tau) \right|^2 $ is integrable
on $[a,b]$, and
so is finite for all $\tau$ in a set of full Lebesgue measure
$P_{\vh}$. Thus the $n$-th term tends to zero for all $\tau\in
P_{\vh}$. Intersecting over all 
integer vectors
$\vh\neq 0$ we get one set $P$ of full measure which works for
all $\vh\neq 0$, that is for all $\tau\in P$
$$
 S_{\vh}(n^2,\tau) \to 0, \qquad n\to \infty \;.
$$

Now given any $N$, find $n$ so that $n^2\leq N<(n+1)^2$. Writing
$N=n^2+k$, $0\leq k\leq 2n$ we have by using the trivial bound
$|e^{2\pi i x}|\leq 1$ that 
$$
\left| S_{\vh}(n^2+k,\tau)-S_{\vh}(n^2,\tau) \right| \ll \frac k{n^2} \ll
\frac 1n \to 0,\qquad n\to \infty
$$
for $\tau\in P$  and since $S_{\vh}(n^2,\tau)\to 0$, we get 
$S_{\vh}(N,\tau)\to 0$ for all $N\to \infty$. 

To show 2.2.A, we square out the sum and directly integrate to
get 
$$
\int_a^b \left| S_{\vh}(N,\tau) \right|^2 d\tau = 
\frac 1{N^2}\sum_{k,l\leq N}
e^{2\pi i\vh \cdot (\va_k-\va_l)}
\int_a^b e^{2\pi i (k-l) \vh \cdot \vom \tau} d\tau \;.
$$
The ``diagonal'' terms $k=l$ give a total contribution of 
$(b-a)/N$ to the
sum, so to prove 2.2.A it suffices to bound the off-diagonal
terms $k\neq l$. 

Since $\vom\cdot \vh\neq 0$ for integer $\vh\neq \vz$ (that was the
assumption on $\vom$), each off-diagonal term contributes 
$$
e^{2\pi i\vh \cdot (\va_k-\va_l)} 
\frac{e^{2\pi i (k-l)\vh \cdot\vom b} - 
e^{2\pi i (k-l)\vh \cdot \vom a}}{2\pi i (k-l)\vh \cdot\vom} \;.
$$
Taking absolute values and summing over all pairs $1\leq k\neq l\leq N$ 
gives a contribution bounded by a constant times 
$$
\frac 1{N^2} \sum_{1\leq k\neq l\leq N} \frac 1{|k-l|}. \;\leqno
(2.2.B)
$$
For fixed $n\neq 0$, the number of solution of $k-l=n$ with $1\leq 
k\neq l\leq N$ is $N-|n|$ if $1\leq |n|\leq N-1$ and zero
otherwise. Thus 2.2.B is given by 
$$
\frac 2{N^2}\sum_{n=1}^{N-1} \frac{N-n}n \ll \frac{\log N}N
$$
which proves 2.2.A and the Theorem. 
\QED 

\medskip
\noindent
{\it 2.3 Cyclic subgroups of $S^1$ are kick unstable.}

\medskip
\noindent
We present here a counter-example to kick stability
constructed by Dima Burago. 
In what follows we identify Kronecker maps
$x \to x + b$ with the corresponding elements $b \in S^1$.
Fix an irrational element
$\omega \in S^1$. The corresponding cyclic subgroup
is uniformly distributed in $S^1$. We claim that
there exists a sequence of kicks $\{\beta_i\} \in (S^1)^{\infty}$
such that for every $\tau \in \N$
the evolution of the kicked system 
is {\bf not} uniformly distributed in $S^1$.
Recall from 2.1 that this evolution 
is given by $ f^{(k)} = \alpha_k + k\tau\omega \mod 1$, where
$\alpha_k = \beta_1 +...+\beta_k$. In order to prove the
claim, choose a function $u: \N \to \N$ such that
the preimage $u^{-1}(k) \subset \N$ of every integer $k \in N$
is a subset of strictly positive density. Put now
$\alpha_k = -u(k)k\omega \mod 1$, 
and $\beta_k = \alpha_k - \alpha_{k-1}$.
Fix $\tau \in \N$ and consider the sequence $\{f^{(k)}\}$.
Since every element of this sequence with $k \in u^{-1}(\tau)$
vanishes, and the set $u^{-1}(k)$ has positive density in $\N$,
we conclude that this sequence is not uniformly distributed
in $S^1$. This completes the proof of the claim.

\vfill\eject

\section*{3. Detecting stable mixing on 
$PSL(2,\RR)/\Ga$}

In this section we prove Theorems 1.5.B,D and in an Appendix, 
present some partial answers on Question 1.8.C.  

\medskip
\noindent
{\it 3.1. Proof of 1.5.B:}

\medskip
\noindent
Let $C$  be the
constant from Definition 1.4.A.  Put $C_1=\max_{j\in\NN} 
\rho (\phi^{-1}_j)$.  The maximum is finite
since $\{\phi_j\}$  represents a finite number
of conjugacy classes, and $\rho$  is bi-invariant
up to $C$  (see 1.4.A). Denote by
$f^{(k)}(\tau)$  the evolution of the kicked
system,
$$f^{(k)}(\tau)=\phi_kh^\tau\ldots\phi_1h^\tau\ .$$
Note that for every $k>0$  and
$\tau\in(0+\infty)$  
$$h^{\tau k}=f^{(k)}(\tau)\cdot\prodl^k_{j=1}h^{-\tau j}
\phi^{-1}_jh^{\tau j}\ .$$
Applying $\rho$  to both sides of this equation
and using properties listed in 1.4.A we get that
$$\rho (h^{\tau k})\le\rho (f^{(k)}(\tau))+\suml^k_{j=1}
\rho(\phi^{-1}_j)+2Ck\leqno(3.1.A)$$
$$\phantom{day}\le \rho
(f^{(k)}(\tau))+k(2C+C_1)\ .$$
Choose $\tau_0 >0$  so large that $\rho (h^{\tau
k})\ge 0.5\tau k\rho_\infty(h)$ for all $\tau
>\tau_0$  and $k\in\NN$.  Put
$C_2=0.5\tau_0\rho_\infty (h)-(2C+C_1)$.
Increasing if necessary $\tau_0$  we assume that
$C_2 >0$.  In view of (3.1.A) we have $\rho
(f^{(k)}(\tau))\ge C_2k$  for    $\tau
>\tau_0$.  Thus for $\tau\ge\tau_0$  the
sequence $f^{(k)}(\tau)$  goes to infinity (see
1.5). Applying the Howe-Moore theorem 1.5.A
we see that for $\tau >\tau_0$  the kicked
system is mixing.\hfill $\Box$

\medskip
\noindent
{\it 3.2.  Quasi-morphisms} 

\medskip
\noindent
Our purpose is to prove Theorem 1.5.D, that is 
if $G\subset PSL(2,\RR)$ is a discrete group, and  $g\in G$ an element
of infinite order, then the following are equivalent:
\begin{itemize} 
\item[(i)] There exists a quasi-morphism $r:G\to \R$ so that
$r_\infty(g)\neq 0$. 
\item[(ii)] $g$ is not conjugate to 
its inverse $g^{-1}$ (in which case we
say that $g$ does not admit a time-reversing symmetry in $G$).
\end{itemize} 

One direction is immediate: Given a quasi-morphism $r$, we note that
$r_\infty$ is {\em homogeneous}: $r_\infty(g^k)=k r_\infty(g)$ for
$k\in \Z$ and so
if $g$ has finite order then clearly $r_\infty(g)=0$. 
Moreover, if $g=hg^{-1} h^{-1}$ with $h\in G$ then also 
$g^k=hg^{-k} h^{-1}$ for all $k\geq 1$. Since $r_\infty$ is a
homogeneous quasi-morphism, we have 
$$
r_\infty(g^k) = r_\infty(hg^{-k} h^{-1}) = r_\infty(g^{-k}) + O(1) = 
-kr_\infty(g) + O(1)
$$
and consequently $2k|r_\infty(g)| = O(1)$ is bounded for all $k\geq
1$, which forces $r_\infty(g) = 0$. 

Let us show now that (ii) yields (i). 
Thus if $G$ is a discrete subgroup of $PSL(2,\RR)$, 
we wish to show that given any element of $G$ of infinite order which
is not conjugate in $G$ to its inverse, there is a homogeneous 
quasi-morphism $r=r_\infty$
of $G$ for which $r_\infty(g)\neq 0$.  

\medskip
\noindent
{\bf Remark.}  
Concerning the condition that $g$ is conjugate in $G$ to
its inverse $g^{-1}$, we note that this can only happen for {\em
hyperbolic} $g$ since elliptic and parabolic  elements of $PSL(2,\RR)$
are never conjugate to their inverses in $PSL(2,\RR)$. 
Moreover, 
if $g\in PSL(2,\RR)$ is hyperbolic then it can be shown that any time
reversing symmetry $K$ of $g$ (i.e. an element $K$ such that $g=Kg^{-1}
K^{-1}$) must satisfy  $K^2=-1$ (in $SL(2,\RR)$). 
Thus for many cases of interest, such as surface groups, which do not
have elliptic elements, this possibility does not arise.

\medskip
\noindent 
We will use the following well-known construction (see \cite{Barge} and 
\cite{Picaud}, 3.3.2). 
For a pair of points $x,y$ in the hyperbolic 
plane $\HH$, write $\ell(x,y)$
for the oriented geodesic 
joining $x$ to $y$. Let $\Omega=dx \wedge dy/y^2$
be the hyperbolic area form. 

\proclaim Definition 3.2.A. 
A one-form $\alpha$ on $\HH$ is {\em bounded} if there is some $C>0$
so that $|\frac {d\alpha}{\Omega}|\leq C$.

Given a bounded $G$-invariant one-form on $\HH$ and a base-point $x\in
\HH$, set 
$$
r_x(g) = r_{x}^{\alpha}(g):= \int_{\ell(x,gx)} \alpha
$$

\proclaim Lemma 3.2.B. 
If $\alpha$ is a bounded $G$-invariant one-form on $\HH$ then 
\begin{itemize}
\item[(i)] $r_x$ is a quasi-morphism.
\item[(ii)] $|r_x-r_y| \leq C_\alpha$. 
\end{itemize}

\noindent{\bf Proof:} 
1) Let $g,h\in G$ and consider 
$$
\delta r_x(g,h) :=r_x(g)+r_x(h)-r_x(gh) \;,
$$
which we want to show is bounded. By $G$-invariance of $\alpha$, we
have 
$$
r_x(h):=\int_{\ell(x,hx)}\alpha = \int_{\ell(gx,ghx)} \alpha
$$
Therefore 
$$
\delta r_x(g,h) =\left(\int_{\ell(x,gx)} + \int_{\ell(gx,ghx)}
-\int_{\ell(x,ghx)}\right) \alpha = \int_{\partial T} \alpha
$$ 
is the integral around the oriented boundary of the geodesic triangle
$T$ with vertices   at $x$, $gx$ and $gh x$. By Stokes' theorem, this
equals the integral of $d\alpha$ on $T$, and thus if $|d\alpha|\leq
C\Om$ then 
$$
|\delta r_x(g,h) |=\left| \int_T d\alpha \right| \leq C\left| \int_T
\Om\right| = C\cdot \mathrm{ area}(T) \;.
$$
Since the area of a geodesic triangle in the hyperbolic plane is at most
$\pi$, we find that $|\delta r_x(g,h) |\leq \pi C$ is bounded and thus
$r_x$ is a quasi-morphism. 

2) To see independence of $r_x$ on the base-point up to a bounded
quantity, consider the integral of $\alpha$ over the boundary of the
geodesic parallelogram $P$ with vertices at $x$, $gx$, $gy$ and $y$, and
again use Stokes' theorem: 
$$
\left( \int_{\ell(x,gx)} +
\int_{\ell(gx,gy)}+\int_{\ell(gy,y)}+\int_{\ell(y,x)} \right) \alpha 
=\int_P d\alpha
= O(1) \;.
$$
By $G$-invariance of $\alpha$, we have 
$$
\int_{\ell(gx,gy)}\alpha  = \int_{\ell(x,y)}\alpha  = 
-\int_{\ell(y,x)}\alpha 
$$
and so we find 
$$
\left| r_x(g) - r_y(g)\right| =
\left| \int_{\ell(x,gx)} +\int_{\ell(gy,y)} \right| = \left| \int_P
d\alpha \right| \leq 2\pi C
$$
since the hyperbolic area of $P$ is at most $2\pi$. 
\QED

Since $G$ is discrete, there are two kinds of elements of $G$ with
infinite order: hyperbolic and parabolic. 
The construction of $r_\infty$ is carried out 
separately for each of these two
cases. 

\medskip
\noindent
{\bf The hyperbolic case.}  
It suffices to consider the case that $g$ is a {\em primitive}
hyperbolic element of $G$, that is we cannot write $g=g_1^k$ for some
$g_1\in G$ and $|k|\geq 2$. Thus we assume this to be the case from
now. 

We recall some facts from the geometry of discrete groups: 

Any hyperbolic element leaves invariant a unique geodesic in $\HH$. 
Let $L$ be the invariant geodesic for a
primitive hyperbolic element $g$.
The following  is a  standard fact: 

\proclaim Lemma 3.2.C. 
\begin{itemize}
\item[(i)] Suppose $\gamma \in G$, $\gamma L=L$ and $\gamma$
preserves the orientation of $L$. Then $\gamma
=g^k$ for some integer $k$. 
\item[(ii)] Suppose  $\gamma \in G$, $\gamma L=L$ and $\gamma$
reverses  the orientation of $L$. Then $\gamma g\gamma^{-1}=g^{-1}$. 
\end{itemize}

Any discrete group $G$ admits a fundamental region $D$ for which the
tessellation $\{\gamma \overline D: \gamma \in G\}$ of the 
upper half-plane $\HH$ is {\em locally finite}, that is to say 
each compact subset of $\HH$ intersects only finitely many of the
translates $\gamma \overline D$.  
An example is the Dirichlet fundamental region
of $G$ \cite{Beardon}.

\proclaim Lemma 3.2.D.  
There is a segment ${\mathcal I}\subset L$ and 
$\epsilon>0$ such that for every
$\gamma\in G$, either $\gamma {\mathcal I}\subset L$ or else 
$dist(\gamma {\mathcal I},L) >\epsilon$. 

\noindent{\bf Proof:}
Let $D$ be a locally finite fundamental domain of $G$ which intersects
$L$ at an interior point. 
Write $U_\delta$ for the $\delta$-neighborhood of $L$ in the
hyperbolic metric (the so-called ``hypercycle domain''). 
Because $D$ is locally finite, 
there are only a {\em finite} number of
distinct translates $\gamma U_\delta$, $\gamma \in G$, which intersect
the closure $\overline D$ of $D$  (see e.g. 
\cite[Theorem 9.2.8~(iii)]{Beardon}). 

Let $Y$ be the union of
these finitely many translates of $U_\delta$, which are
distinct from $U_\delta$ itself. Decreasing
$\delta$, we can guarantee that  $Y\cap L\cap D \neq L\cap D$. 

Choose a segment ${\mathcal I}\subset L\cap D$ such that 
${\mathcal I}\cap Y =  \emptyset$. Assume that for some $\gamma\in G$, 
$\gamma {\mathcal I}\cap U_\delta \neq \emptyset$. 
Then ${\mathcal I}\cap \gamma^{-1} U_\delta \neq \emptyset$. 
But this means that $\gamma^{-1} U_\delta = U_\delta$ due to our
construction. This implies that $\gamma^{-1} L = L$ so that 
$\gamma {\mathcal I}\subset L$. This proves the lemma (with $\epsilon =
\delta$). 
\QED

\noindent{\bf Proof of the theorem in the hyperbolic case: }
Assume that $g$ does not admit a time-reversing symmetry. Then Lemmas 
3.2.C and 3.2.D imply that if $\gamma\neq g^k$ for some $k\in \Z$ then 
$\gamma \mathcal I$ is bounded away from $L$. 
Shrinking if necessary the segment $\mathcal I$, we see that
there exists its small neighborhood
$\mathcal U$ 
such that $\gamma \mathcal U$ is
bounded away from $\mathcal U$ for all $1\neq \gamma \in G$.   

Choose a one-form $\alpha_0$ in $\HH$ such that
$\alpha_0$ has compact support, contained in $\mathcal U$ and
$$\int_{\mathcal I}\alpha_0 >0.$$
Let $\alpha$ be the extension of $\alpha_0$ by periodicity to the
translates $\cup_{\gamma \in G}\gamma \mathcal U$. 
Note that $\alpha$ is a {\em bounded} one-form. 

Now fix a point $x\in L\cap D$ such that
${\mathcal U} \cap L$ lies in the interior of the geodesic
segment $\ell(x,gx)$,
and let
$r=r_x^\alpha$ be the quasi-morphism constructed above.  
Clearly $r(g^n) = \int_{\ell(x,g^nx)} 
\alpha = n\int_{\mathcal I}\alpha_0$ 
and so $r_\infty(g)=\lim r(g^n)/n = \int_{\mathcal I}\alpha_0 >0$. 
This completes the proof. 

\bigskip
\noindent 
{\bf The parabolic case.} 
Let $h\in G$ be a  parabolic element, which as in the hyperbolic case
we may assume is {\em primitive}. Recall that parabolic elements
are never conjugate to their inverses already in $PSL(2,\RR)$.
Let $L$ be an $h$-invariant 
horocycle and let 
${\mathcal U}_L$ be a horocyclic domain in $\HH$, that is a
neighborhood of the cusp fixed by $h$. For instance, if $\infty$ is a
cusp for $G$ and $h(z)=z+1$ and we can take $L=\{y=C\}$ and
${\mathcal U}_L = \{y\geq C\}$. 
The following is well-known: 

\proclaim Lemma 3.2.E. 
One can choose $L$ so that $\gamma {\mathcal U}_L\cap {\mathcal U}_L 
= \emptyset$ for all $\gamma \in G$ with $\gamma \neq h^k$. 

\noindent{\bf Proof:}
We may use a normal form for $h$ and so assume that the cusp is at
$\infty$ and that $h(z)=z+1$. Then if
$\gamma\in G$  does not fix the cusp $\infty$, then 
$\gamma = \left(\smallmatrix a&b\\c&d\endsmallmatrix\right)$ with
$|c|\geq 1$ (see e.g. \cite[Proof of 9.2.8 (ii)]{Beardon}).  
In that case  
$\im(\gamma z) = \im(z)/|cz+d|^2 \leq 1/y$. Thus the horocycles
$L=\{y=C\} $ for $C>1$ satisfy the conditions of the Lemma. 
\QED

\noindent{\bf Proof of the theorem in the parabolic case:}  
For simplicity we assume that $h(z)=z+1$, $z=x+iy$. 
Take $L=\{y=2\}$ and let
$u(y)$ be a smooth cutoff function, with $u(y)\equiv 1$ for $y\geq 3$,
and $u(y)\equiv 0$ if $y\leq 2.5$. Set $\alpha_0=u(y) dx$.

Note that $\alpha_0$ is {\em bounded} on $\HH$, 
since $d\alpha_0 \equiv 0$ for $y\notin (2.5,3)$, while for $2.5<y<3$ we have
$d\alpha_0 = u'(y)dy\wedge dx$ and comparing with 
$\Om = y^{-2}dx\wedge dy$ gives $|d\alpha_0 /\Om|=y^2 u'(y)$ is bounded. 
Moreover $\alpha_0$ is supported in ${\mathcal U}_L=\{y\geq 2\}$. 

Now let $\alpha = \sum_{\gamma\in G/\langle h\rangle}\gamma^*
\alpha_0$ be the periodization of $\alpha_0$. The translates
$\gamma^*\alpha_0$ are supported in distinct translates $\gamma
{\mathcal U}_L$  for distinct $\gamma$ modulo translates by powers of
$h$. Thus we get a $G$-invariant, bounded one-form on $\HH$ 
which equals $dx$ on $\{y \geq 3\}$.

Choose $z$ with $\im(z)=3$ and consider the quasi-morphism
$r=r^\alpha_z$. Then $h^n(z)=z+n$ and clearly 
$\int_{\ell(z,h^nz)} \alpha = n$.  Thus $r_\infty(h)=1>0$ as required.
This completes the proof of Theorem 1.5.D.
\QED

\medskip
\noindent
{\bf Remark 3.2.F.}The phenomenon described in Theorem 1.5.D 
holds true in a more general context of Gromov hyperbolic groups.
In fact, if $G$ is a non-elementary Gromov hyperbolic group
then for every $g \in G$ one has the following alternative.
Either some positive power of $g$ is conjugate to its inverse,
or $G$ admits a homogeneous quasi-morphism which is
positive on $g$. This follows with minor extra efforts
from a work by Epstein and
Fujiwara (see \cite{Epstein}, proof of Lemma 3.5). Another natural
generalization of 1.5.D is as follows. One asks whether
there exists
a quasi- morphism which attains prescribed values on
a given finite subset of the group $G$. A solution of
this problem for discrete subgroups of $PSL(2,\R)$ as well
as an application to stable mixing of linear maps of the 2-torus
will be presented in a forthcoming paper \cite{PR}.

\bigskip
\noindent
\centerline{\bf APPENDIX: The continuous case}

\medskip
\noindent
This appendix to \S 3 is devoted to discussion 
of the continuous case, that 
is we take $G$ to be $PSL(2,\RR)$ acting on $PSL(2,\RR)/\Ga$, where
$\Ga$ is a lattice.  We will take the subgroup  
$$h^t = 
\left (
\begin{array}{cc}
1& t\\
0& 1\\
\end{array}
\right )
$$
which gives the horocycle flow on $X$. Our problem is (cf. 1.8.C
above):

\noindent 
{\it  Is the horocycle flow on
$PSL(2,\RR)/\Ga$  stably mixing?}

\medskip
\noindent
{\it 3.3. Quasi-mixing.}

\medskip
\noindent
First of all, let us
relax the mixing property as follows.

\medskip
\noindent
{\bf Definition 3.3.A.}  A sequential
system $f_*$  acting on a measure space
$(X,\mu)$  by measure-preserving automorphisms
is called {\it quasi-mixing\/} if there exists a
sequence of positive integers $i_k\to +\infty$
such that for any $L^2$-functions $F$  and $H$ on $X$
$$\intl_X F(f^{(i_k)}x)H(x)d\mu\longrightarrow\intl_X
F(x)d\mu\intl_XH(x)d\mu$$
when $k\to\infty$. That is, the subsequence $\{f^{(i_k)}\}$ is mixing.

\medskip
\noindent
Let $\Ga\subset PSL(2,\RR)$  be a lattice.
Consider the left action of $PSL(2,\RR)$ on the
space $PSL(2,\RR)/\Ga$ endowed with the Haar
measure.
Write $h^t=\begin{pmatrix}
1&t\\
0&1
\end{pmatrix}$ for the horocycle flow, and let
$\phi_*=\{\phi_i\}$ be an arbitrary
sequence of kicks from $PSL(2,\RR)$.  Denote by
$QM(\phi_*)$  the set of those periods
$\tau\in (0;+\infty)$  for which the kicked
system $\{\phi_ih^\tau\}$  is quasi-mixing.

 From general considerations, it follows that a quasi-mixing sequence 
is {\em unbounded} (that is has non-compact closure in
$PSL(2,\RR)$). 
It follows from the Howe-Moore theorem 1.5.A that the converse is also true. 

\medskip
\noindent
{\bf Question 3.3.B.}  Is it true that
for every sequence $\phi_*$ the set
$QM(\phi_*)$  has ``large measure"?

\medskip
\noindent
Write $\phi_i=\begin{pmatrix}
a_i&b_i\\
c_i&d_i
\end{pmatrix}$. We give a partial affirmative
answer to Question 3.3.B in terms
of the sequence $\{c_i\}$. 

\medskip
\noindent
\proclaim Theorem 3.3.C.  If $c_i=0$ for all $i$
then the set $(0; +\infty)\bks QM(\phi_*)$     
contains at most $1$ point.

\medskip
\noindent
\proclaim Theorem 3.3.D.  Assume that $c_i\not= 0$
for all $i$  and the sequence ${1\over
n}\suml^n_{i=1}\log |c_i|$ is bounded from
below.  Then the set $(0; +\infty)\bks
QM(\phi_*)$  has finite measure.  

\medskip
\noindent
We prove
these theorems in 3.5, 3.6 below. 

Interestingly enough, Theorems 3.3.C and 3.3.D
handle two opposite cases:  when all $c_i$ vanish,
and when all $c_i$  are bounded away from $0$.
At present it is unclear how to attack the
intermediate situation.

\medskip
\noindent
{\bf Remark 3.3.E.} Let
$M(\phi_*)$  
be the set of those periods $\tau$  for which
the kicked system is mixing (equivalently, goes to
infinity in $PSL(2,\R)$).  In general, one
cannot hope that this set has large measure.
Here is an example which in fact reflects geometry
of real numbers and has nothing to do with the M\"obius
group. We will produce a sequence of kicks of the
form $h^{\beta_i}$, where the sequence $\{\beta_i\}$
is chosen as follows. It is not hard
to exhibit a sequence of intervals
$$I_k = [r_k;r_k + {1\over k}] \subset [0;+\infty), \; k \in
\mathbb{N}$$
which cover every non-negative real number infinitely many
times 
\footnote{
Indeed, one first partitions the divergent series $\sum 1/k$ into
infinitely many divergent subseries $a_{m,n}$, $\sum_n a_{m,n} =
\infty$ for all $m$. To do this, first divide the sequence $1/k$ into
consecutive blocks so that the sum of elements in each block is at
least $1$. 
Then by taking  a bijection $j:\mathbf N\times \mathbf N\to \mathbf N$
one defines  the $m$-th subsequence as the union of the $j(m,n)$-th 
blocks, $n=1,2,\dots $. 

Now denote by $s_{m,n} = \sum_{t\leq n}a_{m,t}$ ($s_{m,0}:=0$)
the partial sums of the $m$-th subsequence, 
and construct the sequence of intervals 
$J_{m,n} = [s_{m,n-1},s_{m,n}]$. Then for each $m\geq 1$, 
$\cup_{n\geq 1} J_{m,n} = [0,\infty)$ and so we get a sequence of intervals
$I_k=[r_k,r_k+1/k]$ which cover every point infinitely many times. }.

Put now $\beta_k = (k-1)r_{k-1} - kr_k$, where
$r_0 = 0$. One calculates that
the evolution of the kicked system is given by
$f^{(k)}(\tau) = h^{k(\tau - r_k)}$. Pick up any positive
real $\tau$. Note that $\tau \in I_k$ if and only if
$k(\tau - r_k) \in [0;1]$, and thus the last inclusion
holds true for infinite number of $k$ due to our choice
of intervals $I_k$.
We conclude
that for every value of $\tau$  the kicked
system is neither mixing, nor its evolution goes
to infinity.  We still have no answer to the
following question.  

\medskip
\noindent
{\bf Question 3.3.F.} Assume that the sequence
of kicks $\phi_* =\{\phi_i\}$ has compact
closure in $PSL(2,\RR)$.  Is it true that the
set $M(\phi_*)$  has large measure? 

\medskip
\noindent
{\it 3.4. A link to discrete Schr\"{o}dinger
equation.}

\medskip
\noindent
Question 3.3.B turns out to be nontrivial even
in the case when the kicks $\phi_i$  have a
very simple form
$$\phi_i=\begin{pmatrix}
1&0\\
c_i&1
\end{pmatrix}
\ ,$$
that is $\phi_i$'s are time-$c_i$-maps of the
conjugate horocycle flow
$$\begin{pmatrix}
0&1\\
-1&0
\end{pmatrix}
h^{-t}
\begin{pmatrix}
0&-1\\
1&0
\end{pmatrix}
\ .
$$
Fix $\tau >0$, and write
$$f^{(k)}(\tau)=\begin{pmatrix}
\al_k&\be_k\\
\ga_k&\de_k
\end{pmatrix}
\ ,
$$
where $f^{(k)}(\tau)$  is the evolution of the
kicked system.  A straightforward calculation
shows that the matrix coefficients satisfy the
following recursive relations:
\begin{eqnarray*}
&&\begin{cases}
\al_k=\al_{k-1}+\tau\ga_{k-1}\\
\ga_k=\ga_{k-1}+c_k\al_k
\end{cases}
\ ,\\[1em]
&&\begin{cases}
\be_k=\be_{k-1}+\tau\de_{k-1}\\
\de_k=\de_{k-1}+c_k\be_k
\end{cases}
\ .
\end{eqnarray*}

Both sequences $\{\al_k\}$ and $\{\be_k\}$
satisfy the second order difference equation
$$q_{k+1}-(2+\tau c_k)q_k+q_{k-1}=0\ \quad k\ge
1\ .\leqno(3.4.A)$$ 
Note that this is the discrete Schr\"{o}dinger
equation with a potential, which depends on the
parameter $\tau$.
Every solution of 3.4.A is uniquely determined
by the initial conditions $q_0$ and $q_1$.  We
get the following result.

\medskip
\noindent
\proclaim Proposition 3.4.B.  The sequence
$\{f^{(k)}(\tau)\}$, $k\in\NN$  is bounded if and
only if all the solutions of the Schr\"{o}dinger
equation (3.4.A) are bounded.

\medskip
\noindent
It is instructive to translate Proposition 3.4.B
into the language of operator theory.  Consider
the space $V$  of all real sequences
$q=(q_1,q_2,\ldots)$.  Define linear operators
on $V$, 
$$L:(q_1,q_2\nek q_i,\ldots)\mapsto (c_1q_1,c_2q_2\nek
c_iq_i,\ldots )$$
and
$$\De_u:(q_1,q_2\nek q_i,\ldots )\mapsto (uq_1
-q_2,-q_1+2q_2-q_3\nek
-q_{i-1}+2q_i-q_{i+1},\ldots )\ .$$
Here $u\in\RR$  is a parameter, and the $i$-th
coordinate of $\De_uq$  is simply the second
difference of the sequence $q$  for all $i\ge 2$. Consider
an operator $K_{u,\tau}=\tau L+\De_u$.  Note
that every vector $q\in Ker K_{u,\tau}$
describes a solution of the Schr\"{o}dinger
equation 3.4.A with the boundary condition
$q_0=(2-u)q_1$. 

Consider the subspace $V_b\subset V$ consisting
of all bounded sequences.  With this notation
the discussion above leads to the following
statement.

\medskip
\noindent
\proclaim Proposition 3.4.C.  The sequence
$\{f^{(k)}(\tau)\}$, $k\in\NN$  is bounded if
and only if $Ker K_{u,\tau}\subset V_b$  for all
$u\in\RR$.

\medskip
\noindent
The next result is a version of Theorem 3.3.D above.

\medskip
\noindent
\proclaim Proposition 3.4.D.  Suppose that
$|c_i|\ge\eps >0$ for all $i\in\NN$.  Then there
exists $\tau_0>0$  such that the sequence $\{
f^{(k)}(\tau)\}$  is unbounded for all $\tau
>\tau_0$.

\medskip
\noindent
{\bf Proof:} Fix $u\in\RR$. We claim that
for $\tau$  large enough $Ker K_{u,\tau}\cap
V_b=\{ 0\}$.  Assume the claim.  Since $\dim Ker
K_{u,\tau}=1$  we get that $Ker K_{u,\tau}$  is
not contained in $V_b$  for $\tau$  large
enough.  Thus the desired result follows from 3.4.C.  

It remains to prove the claim.  Endow the space
$V_b$ with the norm $\| q\| =\sup_i|q_i|$.  Our
assumption on $c_i$  implies that operator $L$
is invertible, $L^{-1}(V_b)=V_b$ and
$\|L^{-1}\|\le {1\over\eps}$.  Further, $\De_u(V_b)
\subset V_b$,  and $\De_u$  is a bounded operator.
Denote $\|\De_u\|=v$.  We have to solve the
equation $K_{u,\tau}q=0$, $q\in V_b$.  Note then
that $(\tau L+\De_u)q=0$,  that is $(\done
+\tau^{-1}L^{-1}\De_u)q=0$.  Since
$\|L^{-1}\De_u\|\le v/\eps$ we see that the
operator $\done +\tau^{-1}L^{-1}\De_u$  is 
invertible for $\tau >v/\eps$, and therefore
$q=0$.  This proves the claim. \hfill $\Box$

\medskip
\noindent
The proof above illustrates the difficulty which
one faces in the case when the coefficients
$c_i$  are allowed to approach arbitrarily close
to $0$.  Indeed, the operator $L^{-1}\De_u$
becomes unbounded, and one loses control on the
kernel of $\done +\tau^{-1}L^{-1}\De_u$  even
for large values of $\tau$.  

\medskip
\noindent
Let us present two additional cases when one
gets the affirmative answer to Question 3.3.B
assuming that $\phi_i=\begin{pmatrix}
1&0\\
c_i&1
\end{pmatrix}$.

\medskip
\noindent
{\bf (3.4.E)} $c_i\to 0$ when $i\to\infty$;

\medskip
\noindent
{\bf (3.4.F)} all $c_i$ are non-negative.

\medskip
\noindent
Indeed assume that $c_i\to 0$.  We write $|\psi
|$  for the Euclidean norm of a matrix $\psi\in
PSL(2,\RR)$, $|\psi |=\sqrt{tr\psi\psi^*}$.  Let
us show that the sequence $|f^{(k)}(\tau)|$  is
unbounded for every $\tau >0$.  Assume on the
contrary that for some $\tau >0$  holds $|f^{(k)}(\tau)|
\le K$  for all $k\in\NN$.  Since $c_i\to 0$
there exists $i,j>0$  such that  
$$|\phi_{i+j}h^\tau\cdots\phi_{i+1}h^\tau|\ge
2K^2\ .$$
But
$$|\phi_{i+j}h^\tau\cdots\phi_{i+1}h^\tau
|=|f^{(i+j)}(\tau)\cdot (f^{(i)}(\tau))^{-1}|\le K^2\ .$$
This contradiction proves the claim.
\footnote{This argument was suggested to us by D.~Kazhdan.}

Assume now that all $c_i$ are non-negative, and
for some $\tau >0$  the sequence
$|f^{(k)}(\tau)|$  is bounded.  Write
$f^{(k)}(\tau)=\begin{pmatrix}
\al_k&\be_k\\
\ga_k&\de_k
\end{pmatrix}$, and note that our assumption
implies that the matrix coefficients are non-negative
and bounded above.  The recursive relations
listed in the beginning of this section show
that the sequence
$\{\al_k\},\{\be_k\},\{\ga_k\},\{\de_k\}$  are
non-decreasing, and thus they converge to some
values
$\al_\infty,\be_\infty,\ga_\infty,\de_\infty$.
Since the matrix $\begin{pmatrix} 
\al_\infty&\be_\infty\\
\ga_\infty&\de_\infty
\end{pmatrix}$ belongs to $PSL(2,\RR)$  we have
that either $\al_\infty\neq 0$  or
$\be_\infty\neq 0$.  Assume without loss of
generality that $\al_\infty\not= 0$.  Since
$c_k={\ga_k-\ga_{k-1}\over\al_k}$  we conclude
that $c_k\to 0$  when $k\to\infty$, and we are
in the case (3.4.E) considered above.  This completes
the analysis of (3.4.E) and (3.4.F).

\medskip
\noindent
{\it 3.5. Proof of Theorem 3.3.C.}

\medskip
\noindent
 Write
$\phi_k=\begin{pmatrix}
a_k&b_k\\
0&a^{-1}_k
\end{pmatrix}$.  Fix $k>0$  and for $m\le k$
denote $\psi_m=\phi_k\cdot\ldots\cdot\phi_m$.
Then the evolution of the kicked systems can be
written as follows:
$$f^{(k)}(\tau)=\phi_kh^\tau\cdot\ldots\cdot
\phi_1h^\tau =\Big(\prodl^{k-1}_{j=0}\psi_{k-j}h^\tau
\psi^{-1}_{k-j}\Big)\cdot\psi_1\ .\leqno(3.5.A)$$
Note that if $\psi =\begin{pmatrix}
u&v\\
0&u^{-1}
\end{pmatrix}$, then
$$\psi h^\tau\psi^{-1}=\begin{pmatrix}
1&\tau u^2\\
0&1
\end{pmatrix}\ .\leqno(3.5.B)$$
The matrix $\psi_1$ has the form $\begin{pmatrix}
a_1\cdot\ldots\cdot a_k&w_k\\
0&a_1\cdot\ldots\cdot a_k
\end{pmatrix}$, where $w_k$  is some real number.
Substituting this to (3.5.A) and using (3.5.B)
we get the following expression for the evolution
of the kicked system:
$$f^{(k)}(\tau)=\begin{pmatrix}
a_1\cdot\ldots\cdot a_k&w_k+\tau z_k
\\
0&a^{-1}_1\cdot\ldots\cdot a^{-1}_k
\end{pmatrix}\ ,\leqno(3.5.C)$$ 
where
$$
z_k=
{\suml^k_{i=1}(a_k\cdot\ldots\cdot a_i)^2\over
a_1\cdot\ldots\cdot a_k}.$$

Now assume that for some $\tau_0$, the sequence
$\{f^{(k)}(\tau)\},k\in\NN$ is bounded. Then
there exist constants $\al >\be >0$  such that
$\be\le |a_1\cdots a_k|\le\al$  for all $k$
(look at the diagonal terms of 3.5.C).
Therefore, for every $k$  and $i$, we have
$|a_k\cdot\ldots\cdot a_i|\ge\be \al^{-1}$ and
thus 
$$|z_k|\ge\be^2\al^{-3}k\
.$$
Note that
$$f^{(k)}(\tau)-f^{(k)}(\tau_0)=\begin{pmatrix}
0&r_k(\tau)\\
0&0
\end{pmatrix}\ ,$$
where $r_k(\tau)=
z_k (\tau -\tau_0)$.
Since $|r_k(\tau)|\ge\be^2\al^{-3}k|\tau
-\tau_0|$, we conclude that the sequence
$f^{(k)}(\tau)$ is unbounded for every
$\tau\not=\tau_0$.  This completes the
proof.\hfill $\Box$

\medskip
\noindent
{\it 3.6. Proof of Theorem 3.3.D.}

\medskip
\noindent
We start with the following

\medskip
\noindent
\proclaim Lemma 3.6.A.  Let $p_k$  be a sequence
of real polynomials of degree $n$  with leading
coefficients $\al_k$.  Suppose that
$|\al_k|\ge\la^k$  for some $\la >0$.  Then the
set
$$Y=\{\tau\in\RR|\ \text{the sequence}\
\{p_k(\tau)\}\ \text{is bounded}\}$$
has finite Lebesgue measure.

\medskip
\noindent
{\bf Proof of 3.3.D:}
Write $\phi_k=\begin{pmatrix}
a_k&b_k\\
c_k&d_k
\end{pmatrix}$, and put $p_k(\tau)=$ trace
$(f^{(k)}(\tau))$, where
$f^{(k)}(\tau)=\phi_kh^\tau\cdot\ldots\cdot
\phi_1h^\tau$.  It is easy to see that
$p_k(\tau)$  is a degree $k$  polynomial with
the leading coefficient $c_1\cdot\ldots\cdot c_k$.
Our assumption on the sequence $\{c_i\}$ implies
that $|c_1\ldots c_k|\ge\la^k$  for some $\la >0$.    
Applying Lemma 3.6.A, we get the statement of
the Theorem.\hfill $\Box$

\medskip
\noindent
{\bf Proof of 3.6.A:}  Take $T>0$, and
write $Y_T=Y\cap [-T;T]$.    The sequence ${1\over k}p_k(\tau)$
converges to $0$ when $k\to\infty$  for every
$\tau\in Y_T$.  Applying the Egorov theorem
\cite{Filter} we get that there exists a subset
$Z_T\subset Y_T$  such that  
$${\rm measure} Z_T \ge
{1\over 2} {\rm measure} Y_T$$  and the sequence
${1\over k}p_k(\tau)$  converges {\it
uniformly\/} on $Z_T$.  In particular, there
exists $k_0 >0$  such that $|p_k(\tau)|\le k$
for all $\tau\in Z_T$ and $k\ge k_0$.  
Setting $\tilp_k(\tau)=\al^{-1}_k p_k(\tau)$, this
implies that 
$$|\tilp_k(\tau)|\le k|\alpha_k|^{-1}\quad\text{for
all}\quad \tau\in Z_T\ .\leqno(3.6.B)$$
Note that $\tilp_k(\tau)$ is a polynomial of
degree $k$  with the leading coefficient $1$.  A
theorem due to Polya (see \cite{Timan}, 2.9.13)
states that for
any measurable subset $Z\subset\RR$ 
$$\max_{\tau\in Z}|\tilp_k(\tau)|\ge 2\cdot\quad
\Big({\text{measure}\ Z\over 4}\Big)^k\leqno(3.6.C)\ .$$
Substituting $Z=Z_T$  and combining with (3.6.B)
we get
$$k|\alpha_k|^{-1}\ge 2\cdot\Big({\text{measure}\ Z_T\over
4}\Big)^k\ge 2\cdot\Big({\text{measure}\ Y_T\over
8}\Big)^k\ .$$
Since $|\alpha_k|\ge\la^k$  we obtain that 
$$k\ge 2\cdot\Big({\text{measure}\ Y_T\over
8}\cdot\la\Big)^k\ .$$
This inequality holds for every $k\ge k_0$,  so
$$\text{measure}\ Y_T\le {8\over\la}\ .$$
Since this is true for every $T>0$  we conclude that
$$\text{measure}\ Y\quad \le{8\over\la}\ .$$
This completes the proof.\hfill $\Box$ 

\vfill\eject
\section*{4.~~Stable super-recurrence in Hamiltonian
dynamics}

In the present section we prove slightly more
general versions of Theorems 1.6.A,B,C.
Recall that these theorems provide a sufficient
condition for kick stable energy conservation law and super-recurrence in
terms of Hofer's geometry of the group of
Hamiltonian diffeomorphisms $\Ham (X,\Om)$.  In
practice, it is easier to perform measurements
not on $\Ham (X,\Om)$ itself but on its
universal cover.  It turns out that these simpler
measurements are powerful enough to detect
stable super-recurrence.

At the very end of the section we prove Proposition
1.7.F stated in the Introduction.

\medskip
\noindent
{\it 4.1. Dynamical preliminaries}

\medskip
\noindent

In this section we describe a link between quasi-integrals (see 1.1.C)
and super-recurrence (see 1.1.D) in the context of 
general sequential systems.  Let
$X$  be a compact topological space endowed with
a Borel probability measure $\mu$. Let $f_*=\{ f_i\}$
be a sequential system which acts on $X$  by
$\mu$-preserving homeomorphisms.  Assume that
$f_*$  is not strictly ergodic (see 1.1.B ).  Then there exists a continuous
function $F$ on $X$ with zero mean such that
$F\not\equiv 0$  and 
$$\limsup_{N\to\infty}\max {1\over N}\suml^{N-1}_{i=0}F\circ
f^{(i)}\ge\al\max F \leqno(4.1.A)$$
for some $\al\in (0;1]$.  In this case we say
that $F$  is an $\al$-{\it quasi-integral\/} of
$f_*$.

For an open
subset $A \subset X$ and $N \in \N$ define the counting function
$\nu_{N,A} : X \to \R$ as follows. For $x \in X$ put $\nu_{N,A}(x)$
to be the the cardinality of the set
$$\{i \in [0;N-1]\; {\big |} \; f^{(i)}(x) \in A\}.$$
Define the quantity
$$R({f_*},A) =
{\lim \sup}_{N \to +\infty} \max_{x \in X} {1 \over N} \nu_{N,A}(x).$$
For every $\epsilon > 0$
there exist arbitrarily long finite pieces $\{x_0=x,...,x_{N-1}\}$
of trajectories of $f_*$ which visit $A$ with the frequency at
least
$R({f_*},A) - \epsilon$. Clearly, $A$ is super-recurrent if $R({f_*},A)$
is strictly bigger than $\mu(A)$.

It turns out that in some situations one can extract 
fairly explicit information on super-recurrent sets
from quasi-integrals.
Let $F$ be an $\al$-quasi-integral of $f_*$.
Put $\ga =\Big|{\min F\over\max F}\Big|$, and
denote by $A_c,c\in (0;1)$  the set
$$\{x\in X\mid F(x)\ge c\max F\}\ .$$

\medskip
\noindent
\proclaim Theorem 4.1.B.  Suppose that $\ga
<{\al^2\over 4-4\al}$. Then for every $c\in (0;1)$
which satisfies
$$\Big| c-{\al\over 2}\Big|<{1\over 2}\sqrt{a^2+4\al\ga
-4\ga}$$
the set $A_c$  is super-recurrent for the system
$f_*$.  Moreover, $R(f_*,A_c)-\mu(A_c)\ge {\al
-c\over 1-c}-{\ga\over c +\ga}$. 

\medskip
\noindent
The proof is based on the following lemmas.

\medskip
\noindent
{\bf Lemma 4.1.C.} {\sl For every $c\in (0;\al)$}
$$R(f_*,A_c)\ge{\al -c\over 1-c}\ .$$

\medskip
\noindent
{\bf Proof.} Choose arbitrary $\eps >0$. Denote
$$I_N={1\over N}\suml^{N-1}_{i=0} F\circ
f^{(i)}\ .$$
There exists an arbitrarily large positive
integer $N$  such that for some $x_0\in X$
$$I_N(x_0)\ge (\al -\eps)\max F\ .$$
Write $\nu_N$  for the counting function
$\nu_{N,A_c}$  defined in the beginning of this section.  Clearly,
$$NI_N(x_0)\leq  c\max F(N-\nu_N(x_0))+\max
F\cdot\nu_N(x_0)\ .$$
Combining this with the previous inequality we
get that
$$\max_x{\nu_N(x)\over N}\ge {\al-c-\eps\over 1-c}\ .$$
Since this holds true for all $\eps > 0$  and
for an infinite sequence of positive values of
$N$  we conclude that
$$R(f_*,A_c)=\limsup_{N\to +\infty}{\nu_N(x)\over N}\ge
{\al -c\over 1-c}\ .\eqno\Box$$

\medskip
\noindent
{\bf Lemma 4.1.D.} {\sl For every $c\in (0;\al)$
$$\mu(A_c)\le {\ga\over c+\ga}\ .$$
}

\medskip
\noindent
{\bf Proof.}  Note that
$$0=\intl_X Fd\mu =\intl_{A_c} Fd\mu +\intl_{X\bks A_c}
Fd\mu\ge c\max F\cdot\mu(A_c)+\min F\cdot
(1-\mu(A_c))\ .$$
Thus
$$\mu (A_c)\le{-\min F\over c\max F-\min
F}={\ga\over c+\ga}\ .\eqno\Box$$

\medskip
\noindent
{\bf Proof of 4.1.B.} The assumptions of the
theorem guarantee that ${\al -c\over
1-c}>{\ga\over c+\ga}$.  Applying 4.1.C and
4.1.D we get that 
$$R(f_*,A_c)-\mu(A_c)\ge{\al -c\over 1-c}-{\ga\over
c+\ga}>0\ .$$
This completes the proof.\hfill $\Box$

\medskip
\noindent
{\it 4.2 Geometric preliminaries}

\MS
\noindent
Let $(X,\Om)$ be a
closed symplectic manifold.  For a smooth path
$q^t,t\in [a;b]$  of Hamiltonian diffeomorphisms
of $(X,\Om)$ set 
$${\rm length}\ (q^t) = \intl^b_a\max_{x\in X}Q(x,t)dt\ ,$$
where $Q$ is the normalized Hamiltonian
generating the path.  Here we use the following
normalization:
$$\int_X Q(x,t)d\mu =0$$
for all $t\in [a,b]$.  Let $(h^t)$  be a one
parameter subgroup of $\Ham (X,\Om)$ generated
by the time-independent Hamiltonian $H$.  Define
a function $\ell_H:[0;+\infty)\to [0;+\infty)$
by
$$\ell_H(s)=\hbox{inf length}\ (q^t)\ ,$$
where the infimum is taken over all paths
$(q^t)$, $t\in [0;s]$ with the following
properties:
\begin{itemize}
\item [$\bullet$] $q_0=\done, q^s=h^s$;
\item [$\bullet$] the paths $q^t$ and $h^t$  are
homotopic through smooth paths with fixed end
points.
\end{itemize}
Clearly,
$$s\max H\ge \ell_H(s)\ge\rho(h^s)\ ,\leqno(4.2.A)$$
where $\rho$  is the positive path of Hofer's
norm defined in 1.6.  It is easy to see that
$\ell_H(s+t)\le\ell_H(s)+\ell_H(t)$  for all
$s,t>0$, so the limit
$$\ell_\infty(H)=\lim_{S\to\infty}{\ell_H(s)\over
s\max H}$$
exists.  This quantity always belongs to the
unit segment $[0;1]$.

In a number of interesting situations one can
find non-trivial lower bounds for $\ell_\infty(H)$
using tools of modern symplectic topology.  Here
we present such a bound (see 4.2.D below) which
was obtained in \cite{Pimrn},\cite{Pbook}.  Recall that
a submanifold  $L\subset(X,\Om)$  is called Lagrangian
if $\dim L={1\over 2}\dim X$, and the symplectic
form $\Om$ vanishes on $TL$.

\medskip
\noindent
{\bf Definition 4.2.B.} Let $L\subset X$  be a
closed Lagrangian submanifold.  We say that $L$
has {\it the Lagrangian intersection property\/}
if $L\cap\phi(L)\not=\emptyset$  for every Hamiltonian
diffeomorphism $\phi\in\Ham (X,\Om)$.

For example, the {\em equator} of the $2$-sphere
(that is a simple closed curve which divides the sphere
into two discs of equal areas) clearly has the  Lagrangian intersection
property.

Consider the cylinder $T^*S^1$  endowed with
coordinates $r\in\RR$  and $t\in S^1$.  The
standard symplectic form on $T^*S^1$  is written
as $dr\wedge dt$.  Denote by $Z$  the zero
section $\{ r=0\}$.  For a symplectic manifold
$(X,\Om)$  consider the topological
stabilization $(X\times T^*S^1,\Om+dr\wedge dt)$.
If $L$  is a closed Lagrangian submanifold of
$X$  then $L\times Z$  is a closed Lagrangian
submanifold of $X\times T^*S^1$.

\medskip
\noindent
{\bf Definition 4.2.C.} Let $L\subset X$  be a
closed Lagrangian submanifold.  We say that $L$
has {\it the stable Lagrangian intersection
property\/} if $L\times Z$  has the Lagrangian
intersection property in $X\times T^*S^1$.

\medskip
\noindent{\bf Remark:}  
It is easily seen that stable Lagrangian intersection property implies 
the  Lagrangian intersection property.

\medskip
\noindent
In many situations, one can detect the stable
Lagrangian intersection property with the help
of the Floer homology.  Let us give two
examples. 
\begin{itemize}
\item [$\bullet$] A closed Lagrangian
submanifold $L\subset X$  with $\pi_2(X,L)=0$
has stable Lagrangian intersection property;
\item [$\bullet$] The equator of the $2$-sphere, 
which trivially has the 
Lagrangian intersection property, can be shown to in fact have 
the {\em stable} Lagrangian intersection property.
\end{itemize}

We refer to \cite{Pimrn},\cite{Pbook} for further
details and references.

\medskip
\noindent
{\bf Theorem 4.2.D.} \cite{Pimrn},\cite{Pbook} {\sl Let $(h^t)$
be a one parameter subgroup of $\Ham (X,\Om)$
generated by a normalized Hamiltonian function $H$.
Assume that there exists a closed
Lagrangian submanifold $L\subset X$  with stable
Lagrangian intersection property such that
$H(x)\ge C>0$ for all $x\in L$.  Then $\ell_H(s)\ge
Cs$  for all $s>0$, and
$\ell_\infty(H)\ge{C\over\max H}$.}   

\medskip
\noindent
{\it 4.3. Detecting stable super-recurrence}

\medskip
\noindent
 Let
$(X,\Om)$  be a closed symplectic manifold.
  Let $(h^t)$  be a one parameter
subgroup of $\Ham (X,\Om)$ generated by a
time-independent normalized Hamiltonian $H$.
Take an arbitrary  bounded sequence
$\phi_*=\{\phi_i\}$  of Hamiltonian
diffeomorphisms of $(X,\Om)$, and consider the
kicked system $f^\tau_*=\{\phi_ih^\tau\}$.  Put 
$$A_c=\{x\in X\mid H(x)>c\max H\}\ ,$$
where $c\in (0;1)$.

\medskip
\noindent
\proclaim Theorem 4.3.A.  Suppose that
$\ell_H(s)=s\max H$  for all $s>0$.  Then for
every $c,\alpha\in (0;1)$  and $\eps >0$ there exists a
subset $P\subset (0;\infty)$ with the following
properties:  
\begin{itemize}
\item [$\bullet$] the set $(0; +\infty)\bks P$
has finite Lebesgue measure;
\item [$\bullet$]  for every $\tau \in P$ the Hamiltonian $H$
is an $\alpha$-quasi-integral of the kicked system $f^\tau_*$; 
\item [$\bullet$] for every $\tau\in P$  the set
$A_c$  is super-recurrent for the kicked system
$f^\tau_*$  with $R(f^\tau_*,A_c)>1-\eps$.
\end{itemize}

\medskip
\noindent
{\bf Proof of Theorem 1.2.D:}
Since the maximum set of $H$ contains an equator, and
the equator has the stable Lagrangian intersection property
(see 4.2 above) we conclude from 4.2.D and 4.2.A that
$\ell_H (s) = s\max H$ for all $s$. The Theorem follows now
from 4.3.A.
\QED

\medskip
\noindent
The condition $\ell_H(s)=s\max H$  is very
restrictive.  However our technique enables us
to detect a weaker version of kick stable
super-recurrence in a more general situation.
Recall that the {\it density\/} of a subset
$P\subset (0;+\infty)$  is
$$\liminf_{T\to +\infty}{1\over T}\ \text{measure}\ 
(P\cap (0; T])\ .$$
Suppose that $\ell_{\infty}(H) >0$ and choose any $\alpha \in (0;\ell_{\infty}(H))$.
Put $$\tet={\ell_\infty(H)-\al\over 1-\al}.\leqno(4.3.B) $$

\medskip
\noindent
\proclaim Theorem 4.3.C. There exists a subset  $P \subset (0;+\infty)$
of density at least $\theta$
such that for every $\tau \in P$ the Hamiltonian $H$ is an
$\al$-quasi-integral of the kicked system $f^\tau_*$.

\medskip
\noindent
Theorems 4.3.A and 4.3.C are proved in 4.5 below.

\medskip
\noindent
{\bf Proof of Theorem 1.6.A:} 
The Theorem  follows from 4.3.C. Indeed, $\theta \to \ell_{\infty}(H)$
when $\alpha \to 0$. 
\QED

\medskip
\noindent
Let us describe an application of 
Theorem 4.3.C to stable super-recurrence.
Set 
$\ga =\Bigg|{\min H\over\max
H}\Bigg|$. We assume that the following inequality holds:
$$\ga
(4-4\ell_\infty (H)) 
 <\ell_\infty
(H)^2 
\ .\leqno(4.3.D)$$ 
Denote by $P_{c,\de}$  (where $c,\de\in (0;1)$)
the set of those values of the period
$\tau\in(0;+\infty)$  such that the set $A_c$
is super-recurrent for the kicked system
$f^\tau_*$  with $R(f^\tau_*,A_c)-\mu(A_c)>\de$.
Choose arbitrary $\alpha \in (0;\ell_{\infty}(H))$ such that
$$\ga < {{\al^2} \over {4-4\al}}.$$
Pick any $c\in (0;1)$  which satisfies
$$\Big|c-{\al\over 2}\Big| <{1\over 2}
\sqrt{\al^2+4\al\ga -4\ga}\ .\leqno(4.3.E)$$
Set
$$\de ={\al -c\over 1-c}-{\ga\over c+\ga}\leqno(4.3.F)$$
One concludes from (4.3.E) that $\de$  is positive.
The next result 
follows immediately from 4.3.C and 4.1.B.
Here we assume (4.3.D), and define $\theta$ by (4.3.B).

\medskip
\noindent
\proclaim Theorem 4.3.G. The density of $P_{c,\de}$ is 
greater or equal to $\theta$.  

\medskip
\noindent
{\bf Proof of Theorem 1.6.B.} It follows from 4.2.A that
$\ell_{\infty}(H) = 1$. Thus assumption 4.3.D holds.
Put $c=\epsilon$ and choose $\alpha$ so close
to $1$ that 4.3.E holds.
Note that $\theta = 1$.
 Define $\delta$ by 4.3.F. It follows
from 4.3.G that $A_c$ is super-recurrent when $\tau$ belongs
to the set $P_{c,\de}$ of density 1.
\QED

\medskip
\noindent
{\bf Proof of Theorem 1.6.C.} It follows from
(4.2.A) that $\ell_\infty(H)\ge 0.9$.  Further,
it is given that $\ga =1$. 
Thus assumption 4.3.D holds.
 Choose $\al =0.83$,
and $c=0.4$.  Note that the inequality (4.3.E)
is satisfied and $\tet={0.9-0.83\over 1-0.83} >0.4$.
Theorem follows now from 4.3.G.\hfill $\Box$

\medskip
\noindent
{\it 4.4. A geometric inequality}

\medskip
\noindent
The main ingredient of our approach to Theorems
4.3.A and 4.3.C is the following upper bound on
the function $\ell_H(s)$. Let $(h^t)$  be a one
parameter subgroup of $\Ham (X,\Om)$  generated
by a time-independent Hamiltonian $H$ with zero mean.
Take an arbitrary sequence $\{\phi_i\}$  of
Hamiltonian diffeomorphisms.  Denote by $f^{(i)}(\tau)$
the evolution of the kicked system, that is
$$f^{(i)}(\tau)=\phi_ih^\tau\phi_{i-1}h^\tau
\cdot\ldots\cdot\phi_1h^\tau\ ,$$
where $i\ge 1$, $\tau\in(0;\infty)$. Put $f^{(0)}(\tau)
\equiv\done$.  We write $\bar \rho$  for 
Hofer's norm defined 1.2.C. 

\medskip
\noindent
\proclaim Theorem 4.4.A.  For every $N\in\NN$
and $T>0$  holds
$$\ell_H(NT)\le\intl^T_0\maxl_X\sum^{N-1}_{i=0}H\circ
f^{(i)}(t)dt+ 2\sum^{N-1}_{i=1}{\bar \rho}(\phi_i)
\ .$$

\medskip
\noindent
{\bf Proof:} The proof is divided into several
steps.

\medskip
1)  Decompose
$h^{NT}=A_{N,T}B_{N,T}$, where
$$A_{N,T}=h^T\cdot\prodl^{N-1}_{i=1}\psi_ih^T\psi_i^{-1}\
,$$
with 
$$\psi_i=\psi_{i,N} =\phi_{N-1}\cdot\ldots\cdot\phi_{N-i}\
,\quad i=1\nek N-1\ $$
and 
$$B_{N,T}=\phi_{N-1}\cdot\ldots\cdot\phi_1
\prodl^{N-1}_{j=1}h^{-jT}\phi^{-1}_jh^{jT}\ .$$
Take $\eps >0$, and choose paths $\phi^{(s)}_i$,
$s\in [0;T]$ of Hamiltonian diffeomorphisms
which join $\done$  with $\phi_i$ so that the lengths
of the paths
$\phi^{(s)}_i$ and $(\phi^{(s)}_i)^{-1}$ do not exceed 
${\bar \rho} (\phi_i)+\eps$ for all $i$.
Consider the following paths of Hamiltonian
diffeomorphisms defined for $s\in [0;T]$:
\begin{eqnarray*}
&&c_s=h^{Ns}\\
&&a_s=h^s\prodl^{N-1}_{i=1}\psi_ih^s\psi^{-1}_i
\end{eqnarray*}
and 
$$b_s=\phi^{(s)}_{N-1}\cdots\phi^{(s)}_1
\prodl^{N-1}_{j=1} h^{-jT}(\phi^{(s)}_j)^{-1}
h^{jT}\ .$$
The paths $\{ c_s\}$ and $\{ a_sb_s\}$ join $\done$
with $h^{NT}$. 
 
\medskip
2) We claim that the paths $\{c_s\}$  and
$\{a_sb_s\}$  are homotopic with fixed
endpoints.  Indeed take the parameter of
homotopy $u\in [0;1]$  and write
$$a_{s,u}=h^s\prodl^{N-1}_{i=1}\psi^{(u)}_ih^s 
(\psi^{(u)}_i)^{-1}\ ,$$
where $\psi^{(u)}_i=\phi_{N-1}^{(Tu)}\cdots
\phi^{(Tu)}_{N-i}$.  Set
$$b_{s,u}=\phi^{(su)}_{N-1}\cdots\phi_1^{(su)}
\prodl^{N-1}_{j=1}h^{-jT}(\phi_j^{(su)})^{-1}h^{jT}\
.$$
The required homotopy is given by
$d_{s,u}=a_{s,u}b_{s,u}$  where $s\in [0;T]$
and $u\in [0;1]$.

\medskip
3) It follows from the definition of the
function $\ell_H$  that $\ell_H(NT)\le \length
(a_sb_s)_{s\in [0;T]}$.
But length $(a_sb_s)\le\length (a_s)+\length
(b_s)$, and $$\length (b_s)\le 2\suml^{N-1}_{j=1}{\bar \rho}
(\phi_j)+2(N-1)\eps 
$$ 
(one uses here that Hofer's norm $\bar\rho$ is bi-invariant). 
We conclude that 
\begin{itemize}
\item[(4.4.B)] $\qquad\qquad\qquad \ell_H(NT)\le\length
(a_s)_{s\in[0;T]} +$
\item[] $\qquad\qquad\qquad\qquad\quad
+2\suml^{N-1}_{j=1}{\bar \rho}(\phi_j)
+2(N-1)\eps$.
\end{itemize}

\medskip
4) Denote by $\tilF(x,s)=\tilF_s(x)$  the
normalized Hamiltonian function generating
$(a_s)$.  In order to calculate $\tilF$, 
we use the following product
formula:  Let $(p_s),(q_s), s\in [0;T]$  be two
Hamiltonian flows generated by normalized Hamiltonians
$P(x,s)$  and $Q(x,s)$.  Then the product
$(p_sq_s)$  is a Hamiltonian flow generated by
normalized Hamiltonian $P(x,s)+Q(p^{-1}_sx,s)$.
In particular, for a given Hamiltonian diffeomorphism
$\psi$ the path 
$(\psi
h^s\psi^{-1})$ is generated by Hamiltonian $H\circ \psi^{-1}$.

Applying the product formula we get that 
$$\tilF_s=H+\suml^{N-1}_{i=1}H\circ\psi^{-1}_i\circ
\bl\prodl^{i-1}_{j=0}\psi_jh^s\psi^{-1}_j\br^{-1}\ ,$$
where $\psi_0=\done$.  Using that
$\psi^{-1}_{k-1}\psi_k=\phi_{N-k}$ we can
rewrite this as
$$\tilF_s=H+\suml^{N-1}_{i=1}H\circ\bl\prodl^i_{j=1}h^s
\phi_{N-j}\br^{-1}\ .$$
Introduce a new function 
$$F_s=\tilF_s\circ \bl \prodl^{N-1}_{j=1} h^s\phi_{N-j}
\br\circ h^s\ .$$
Then
\begin{eqnarray*}
&&F_s=H\circ h^s\circ \bl \prodl^{N-1}_{j=1} \phi_{N-j}
h^s\br +\\
&&\qquad+\suml^{N-2}_{i=1} H\circ \bl \prodl^{N-1}_{j=i+1}
h^s \phi_{N-j}\br\circ h^s+H\circ h^s\ .
\end{eqnarray*}
The energy conservation law implies that $H\circ
h^s=H$.  Thus
\begin{eqnarray*}
&&F_s=H\circ\prodl^{N-1}_{j=1}\phi_{N-j}h^s+\\
&&\qquad +\suml^{N-2}_{i=1}H\circ\phi_{N-i-1}h_s
\phi_{N-i-2}h_s\cdots\phi_1h_s+H\ .
\end{eqnarray*}
The last expression can be rewritten in terms of
the kicked system:
$$F_s=\suml^{N-1}_{k=0}H\circ f^{(k)}(s)\ .$$
Since $\maxl_X F_s=\maxl_X\tilF_s$  for all $s$,
we get that
\begin{eqnarray*}
&&\length(a_s)_{s\in
[0;T]}=\intl^T_0\maxl_X\tilF_sds=\\
&&\qquad =\intl^T_0\maxl_X\suml^{N-1}_{k=0}H\circ
f^{(k)}(s)ds\ .
\end{eqnarray*}
Substituting this expression into (4.4.B) we get
the desired inequality 4.4.A.  This completes the proof.

\medskip
\noindent
{\it 4.5 Proof of main theorems} 

\medskip
\noindent
{\bf Proof of Theorem 4.3.C.}
Set $H_N(t)={1\over N}\maxl_X\suml^{N-1}_{i=0}
H\circ f^{(i)}(t)$.  Let $P$  be the set of
those $t\in (0;\infty)$  for which the
inequality $H_N(t)\ge\al\max H$  holds for an
infinite sequence of positive integers $N$.  We
have to show that the density of $P$  is at least
$\tet={\ell_\infty (H)-\al\over 1-\al}$.

Since the sequence $\{\phi_i\}$ is bounded, 
there exists $u>0$
such that $2{\bar \rho}(\phi_i) \le u$
for all $i$.  
Fix $\ka >0$, and abbreviate $\ell =\ell_\infty(H)$.
There exist $T_0>0$, $N_0>0$ such that for all
$N>N_0$, $T>T_0$
$$\ell_H(NT)\ge NT(\ell -\ka )\max H\ . $$
The geometric inequality 4.4.A yields
$$\ell_H(NT)\le N\intl^T_0H_N(t)dt+(N-1)u\ .$$
Therefore
$$NT(\ell -\ka)\max H\le N\intl^T_0H_N(t)dt
+Nu\ .$$
Increasing $T_0$  we can assume that $u\le\ka
T\max H$, so we get
$$\intl^T_0H_N(t)dt\ge (\ell -2\ka)T\max
H\ .\leqno(4.5.A)$$
Consider the set
$$Q_{N,T}=\{t\in [0; T]\mid H_N(t) <\al \cdot
\max H\}\ .$$
Clearly,
\begin{eqnarray*}
&&\intl^T_0H_N(t)dt=\intl_{Q_{N,T}}H_N(t)dt+\intl_{[0;T]
\bks Q_{N,T}}H_N(t)dt\le\\
&&\qquad\quad\le\al\max H\cdot\text{measure}(Q_{N,T})+\max
H\cdot (T-\text{measure}(Q_{N,T}))\ .
\end{eqnarray*}
Combining this with 4.5.A we get that
$${\text{measure}(Q_{N,T})\over T}\le {1-\ell
+2\ka\over 1-\al}\leqno(4.5.B)$$

Clearly,
$[0;T]\bks P=\bigcup^\infty_{k=1}\bigcap^\infty_{N=k}
Q_{N,T}$, so in view of (4.5.B)

$${1\over T}\text{measure}(P\cap
[0;T]) \ge{\ell -\al-2\ka\over 1-\al}$$ for
$T>T_0$.  Thus  
$\text{density}(P)\ge\tet ={\ell-\al-2\ka\over 1-\al}$.
Since $\ka$ can be chosen arbitrary small, 
this completes the proof of the
theorem.\hfill $\Box$

\medskip
\noindent
{\bf Proof of Theorem 4.3.A}
The proof is analogous to the one of 4.3.C.
Take $\al\in (0;1)$  and define $H_N(t)$, $u$ and
$Q_{N,T}$ as above.  Since
$\ell_H(NT)=NT\max H$  due to our assumption, and
$$\ell_H(NT)\le N\intl^T_0H_N(t)dt+Nu$$
in view of 4.4.A, we get that 
$$\intl^T_0H_N(t)dt\ge T\max H-u\ .$$
Therefore
$$\text{measure}(Q_{N,T})\le {u\over (1-\al)\max
H}\ .$$
The set $Q=\bigcup_T\bigcup^\infty_{k=1}
\bigcap^\infty_{N=k} Q_{N,T}$ has finite measure.
Consider its complement $P=(0; +\infty)\bks Q$.
For every $\tau\in P$  function $H$  is an
$\al$-quasi-integral of the kicked system
$f^\tau_*$.

Assume without loss of generality that $\alpha$ is sufficiently 
close to 1 so that the following hold:
\begin{itemize}
\item [$\bullet$] $\Big|c-{\al\over 2}\Big|
<{1\over 2}\sqrt{\al^2+4\al\ga -4\ga}$;
\item [$\bullet$] ${\al -c\over 1-c}\ge 1-\eps\
,\ {\al -c\over 1-c} >{\ga\over c+\ga}$.
\end{itemize}

Here $\ga$ stands for $\Big|{\min H\over\max H}\Big|$,
and $\eps$  is given in the formulation of
4.3.A.  It follows from 4.1.B that the set
$A_c=\{x\in X\mid H(x)\ge c\max H\}$  is
super-recurrent for $f^\tau_*$.  Moreover, Lemma
4.1.C guarantees that $R(f^\tau_*,A_c)\ge
1-\eps$.  This completes proof.\hfill $\Box$ 

\medskip
\noindent
{\it 4.6 Creating random behaviour} 

\medskip
\noindent
Here we prove Proposition 1.7.F of the
introduction.  Let $(X,\Om)$  be a closed
symplectic manifold such that $\Ham (X,\Om)$  is
$C^\infty$-closed in $\Diff (X)$.  Note 
that under this assumption there
exists a bounded strictly ergodic sequence
$\psi_*=\{\psi_i\}$  of Hamiltonian
diffeomorphisms.  This is an immediate
consequence of \cite{Pjems}, Th.~1.2.A.    

Let $\tet$ be a time reversing symmetry for $(h^t)$, that
is $\tet h^t\tet^{-1}=h^{-t}$  for all $t\in\RR$.  Define
a sequence of kicks $\{\phi_i\}$  as follows.  Set
$\phi_i=\tet^{-1}$  when $i$ is odd, and
$\phi_i=\psi_k\tet$ when $i$  is even and equals $2k$. 
We claim that
for every $\tau >0$  the kicked system
$f^\tau_*=\{\phi_ih^\tau\}$  is strictly ergodic. 
Indeed, note that $f^{(i)}(\tau)=\phi_ih^\tau\cdot
\ldots\phi_1h^\tau$  equals $\psi^{(k)}$  for
$i=2k$, and $\tet^{-1}h^\tau\psi^{(k)}$ for
$i=2k+1$.  For every continuous function $F$  on
$X$  set $I_N=\suml^{N-1}_{i=0}F\circ f^{(i)}(\tau)$.
We see that for $N=2k$, $I_N =\suml^{k-1}_{i=0}F\circ
\psi^{(i)}+\suml^{k-1}_{i=0}(F\circ\tet^{-1}h^\tau)
\circ\psi^{(i)}$.  Thus the uniform limit
$\lim\limits_{N\to\infty}{1\over N}I_N$  exists
and equals ${1\over 2}\Big(\intl_XFd\mu
+\intl_XF\circ\tet^{-1}h^\tau d\mu\Big)=\intl_XFd\mu$.
Therefore $f^\tau_*$ is strictly ergodic for
every $\tau$. This completes the proof. \hfill
$\Box$

\medskip\noindent
{\bf Example 4.6.A}  We conclude this section
with an example in the spirit of 1.7.G and
1.7.H.  It shows that the phenomenon presented
in Theorem 4.3.A  is a purely Hamiltonian one and
may disappear when one allows kicks which are
symplectic but not necessarily Hamiltonian
diffeomorphisms of $(X,\Om)$.  Consider the $2$-torus
$\TT^2=\RR^2(x,y)/\ZZ^2$ endowed with the symplectic 
form $\Om = dx\wedge dy$.  Consider a one
parameter subgroup $(h^t)$ of Hamiltonian
diffeomorphisms given by
$$h^t(x,y)=(x,y-t\sin 2\pi x)\ .$$
It is generated by a normalized Hamiltonian
function $H(x,y) ={1\over 2\pi}\cos 2\pi x$,
which attains its maximal value on a
non-contractible curve $L=\{ x=0\}$.  Since $L$
has stable Lagrangian intersection property (see
4.2.C above) it follows from Theorem 4.2.D that
$\ell_H(s)\equiv s\max H$.  Thus Theorem 4.3.A
implies stable super-recurrence. On the other
hand the shift $\tet:(x,y)\mapsto (x+{1\over
2},y)$  is a time reversing symmetry of $(h^t)$.
Note that $\tet$ is a symplectic, but not a
Hamiltonian diffeomorphism.  Our proof above
shows that there exists a sequence of
symplectomorphisms $\{\phi_i\}$ such that the
kicked system $\{\phi_ih^\tau\}$ is strictly
ergodic for each value of $\tau$. 

\vfill\eject

\vfill\eject
 
\end{document}